\documentclass[a4paper, 11pt, twoside, titlepage]{article}
\usepackage{vmargin}
\usepackage[latin1]{inputenc}
\usepackage{hyperref}
\usepackage{amssymb}
\usepackage{amsmath}
\usepackage{xspace}
\usepackage{oldgerm}
\usepackage{epsfig}
\usepackage{bbm}
\usepackage{makeidx}
\usepackage{theorem}
\usepackage{mathrsfs}
\usepackage{fancyhdr}
\usepackage{xypic}
\usepackage{psfrag}
\usepackage{color}
\usepackage{stmaryrd}
\usepackage{cancel}

\newcommand{\Sym}{\mathrm{Sym}} % Symmetrische Formen
\newcommand{\Ra}{\rightarrow}

\renewcommand{\emph}{\textbf}

\renewcommand{\rho}{\varrho}
\renewcommand{\phi}{\varphi}
\renewcommand{\epsilon}{\varepsilon}

\newcommand{\tr}{\mathrm{trace}}
\newcommand{\Pf}{\mathrm{Pf}} % Pfaffian
\newcommand{\dis}{\mathrm{dis}} % discriminant

\newcommand{\CC}{\mathbbm{C}} % komplexe
\newcommand{\NN}{\mathbbm{N}} % natürliche
\newcommand{\mat}{\mathrm{Mat}} % Algebra der Matrizen
\newcommand{\GL}{\mathrm{GL}} % General Linear Group
\newcommand{\G}{\mathrm{G}} % e.g. G2
\newcommand{\F}{\mathrm{F}} % e.g. F4
\newcommand{\E}{\mathrm{E}} % e.g. E6, E7, E8

\newcommand{\SL}{\mathrm{SL}} % \'Etale Linear Group
\newcommand{\Un}{\mathrm{Un}} % Unipotent Group
\newcommand{\SO}{\mathrm{SO}} % \'Etale Orthogonal Group
\newcommand{\Sp}{\mathrm{Sp}} % symplektische Gruppe
\newcommand{\Spin}{\mathrm{Spin}} % Spin-Gruppe
\newcommand{\lag}{\langle}
\newcommand{\rag}{\rangle}

{\theorembodyfont{\rmfamily} \newtheorem{DEF}{Definition}[section]}
\newcommand{\BDEF}{\begin{DEF}}
\newcommand{\EDEF}{\end{DEF}}
{\theorembodyfont{\rmfamily} \newtheorem{REM}[DEF]{Remark}}
\newcommand{\BREM}{\begin{REM}}
\newcommand{\EREM}{\end{REM}}
{\theorembodyfont{\slshape} \newtheorem{PROP}[DEF]{Proposition}}
\newcommand{\BPROP}{\begin{PROP}}
\newcommand{\EPROP}{\end{PROP}}
{\theorembodyfont{\rmfamily} \newtheorem{EXP}[DEF]{Example}}
\newcommand{\BEXP}{\begin{EXP}}
\newcommand{\EEXP}{\end{EXP}}
{\theorembodyfont{\slshape} \newtheorem{LEMMA}[DEF]{Lemma}}
\newcommand{\BLEMMA}{\begin{LEMMA}}
\newcommand{\ELEMMA}{\end{LEMMA}}
{\theorembodyfont{\slshape} \newtheorem{COR}[DEF]{Corollary}}
\newcommand{\BCOR}{\begin{COR}}
\newcommand{\ECOR}{\end{COR}}
{\theorembodyfont{\slshape} \newtheorem{THM}[DEF]{Theorem}
\newcommand{\BTHM}{\begin{THM}}
\newcommand{\ETHM}{\end{THM}}

{\theorembodyfont{\rmfamily} }

\setcounter{tocdepth}{2}

\makeindex

\begin{document}

\thispagestyle{empty}%

% \maketitle 
 \hspace{1cm} \bigskip 
 
\begin{center}
\LARGE Tables of prehomogeneous modules and\\ \'etale modules of reductive algebraic groups
\end{center}

 \hspace{3cm} \bigskip 
\begin{center}
Compiled by Wolfgang Globke,

based on the results cited in the reference section.
\end{center}
 \hspace{3cm} \bigskip 

%\newpage
\tableofcontents

% ======================
% ======================
\section{Classification of prehomogeneous modules}\label{sec_class}

\subsection{Irreducible reduced prehomogeneous modules}

The irreducible and reduced prehomogeneous modules were classified
by Sato and Kimura, thus we will label each class by
SK $n$, where $n$ is the number given to the class in \S 7 of the
the original work by Sato and Kimura \cite{SK}.
Along with each module, we will state the connected component of the
generic isotropy group, denoted by $G_v^\circ$,
and in some cases the irreducible relative invariant, denoted by $f$.
\index{prehomogeneous module!irreducible, reduced}
\index{prehomogeneous module!classification}
\index{SK}

Let $G$ be a reductive group and $(G,\rho,V)$ an irreducible and
reduced prehomogeneous module. Then it is  equivalent to one
of the following prehomogeneous modules:
\begin{itemize}
% SK I
\item[{\rm \textbf{SK I}}] Regular irreducible reduced prehomogeneous
modules.
\begin{enumerate}
% 1
\item $(G\times \GL_m, \rho\otimes\omega_1, V^m \otimes \CC^m)$,\\
where $\rho:G\Ra\GL(V^m)$ is an $m$-dimensional irreducible
representation of a connected semisimple algebraic group $G$
(or $G=\{1\}$ and $m=1$).
We have $G_v^\circ \cong G$ and
$f(x) = \det(x)$ for $x\in \mat_m\cong V^m\otimes\CC^m$, $\deg(f)=m$.
% 2
\item $(\GL_n,2\omega_1,\Sym^2 \CC^n)$ for $n\geq 2$.\\
We have $G_v^\circ\cong \SO_n$ and $f(x)=\det(x)$ for
$x\in\{A\in\mat_n\mid A^\top = A\}\cong \Sym^2\CC^n$, $\deg(f)=n$.
% 3
\item $(\GL_{2n},\omega_2,\bigwedge^2 \CC^{2n})$ for $n\geq 3$.\\
We have $G_v^\circ\cong\Sp_{n}$ and $f(x)=\Pf(x)$ for
$x\in\{A\in\mat_{2n}\mid A^\top = -A\}\cong\bigwedge^2\CC^{2n}$,
$\deg(f)=n$.
% 4
\item $(\GL_2,3\omega_1,\Sym^3\CC^2)$.\\
We have $G_v^\circ\cong\{1\}$ and
$f(a)=a_2^2 a_3^2 + 18 a_1 a_2 a_3 a_4 - 4 a_1 a_3^3 - 4 a_2^3 a_4
-27 a_1^2 a_4^2$ for
$a=a_1 x^3 + a_2 x^2 y + a_3 x y^2 + a_4 y^3 \in \Sym^3\CC^2$
(so $f$ is the discriminant of a binary cubic form $a(x,y)$).
% 5
\item $(\GL_6,\omega_3,\bigwedge^3 \CC^6)$.\\
We have $G_v^\circ\cong\SL_3\times\SL_3$ and
$f(x)=(x_0 y_0 - \tr(XY))^2 + 4x_0\det(Y)+4y_0\det(X)
-4\sum_{i,j} \det(X_{ij})\det(Y_{ji})$ (see \S 5, p. 83 in \cite{SK}
for a definition), $\deg(f)=4$.
% 6
\item $(\GL_7,\omega_3,\bigwedge^3 \CC^7)$.\\
We have $G_v^\circ\cong \G_2$ and $\deg(f)=7$.
% 7
\item $(\GL_8,\omega_3,\bigwedge^3\CC^2)$.\\
We have $G_v^\circ\cong\SL_3$ and $\deg(f)=16$.
% 8
\item $(\SL_3\times\GL_2,2\omega_1\otimes\omega_1,
\Sym^2\CC^3\otimes\CC^2)$.\\
We have $G_v^\circ\cong\{1\}$ and $f(A,B)=\dis(\det(xA+yB))$ for
$(A,B) \in\{(X,Y)\mid X,Y\in \mat_3, X^\top=X, Y^\top=Y\}
\cong\Sym^2\CC^3\otimes\CC^2$, $\deg(f)=12$.
% 9
\item $(\SL_6\times\GL_2, \omega_2\otimes\omega_1,
\bigwedge^2\CC^6\otimes\CC^2)$.\\
We have $G_v^\circ\cong \SL_2\times\SL_2\times\SL_2$ and
$f(A,B)=\dis(\Pf(xA+yB))$ for
$(A,B)\in\{(X,Y)\mid X,Y\in \mat_6, X^\top=-X, Y^\top=-Y\}
\cong\bigwedge^2\CC^6\otimes\CC^2$, $\deg(f)=12$.
% 10
\item $(\SL_5\times\GL_3,\omega_2\otimes\omega_1,
\bigwedge^2\CC^5\otimes\CC^3)$.\\
We have $G_v^\circ\cong\SL_2$ and $\deg(f)=15$.
% 11
\item $(\SL_5\times\GL_4,\omega_2\otimes\omega_1,
\bigwedge^2\CC^5\otimes\CC^4)$.\\
We have $G_v^\circ\cong\{1\}$ and $\deg(f)=40$.
% 12
\item $(\SL_3\times\SL_3\times\GL_2,\omega_1\otimes\omega_1\otimes\omega_1, \CC^3\otimes\CC^3\otimes\CC^2)$.\\
We have $G_v^\circ\cong\GL_1\times\GL_1$ and $f(A,B)=\dis(\det(xA+yB))$
for $(A,B)\in\mat_3\oplus\mat_3\cong\CC^3\otimes\CC^3\otimes\CC^2$,
$\deg(f)=12$.
% 13
\item $(\Sp_n\times\GL_{2m},\omega_1\otimes\omega_1,
\CC^{2n}\otimes\CC^{2m})$ for $n\geq 2m\geq 2$.\\
We have $G_v^\circ\cong\Sp_{m}\times\Sp_{n-m}$ and
$f(X)=\Pf(X^\top J X)$ for $X\in\mat_{2n,2m}$, $\deg(f)=2m$.
% 14
\item $(\GL_1\times\Sp_3,\mu\otimes\omega_3,
\CC\otimes V^{14})$.\\
We have $G_v^\circ\cong\SL_3$ and $\deg(f)=4$, where $f$ is given by the
restriction of the relative invariant of SK I-5.
% 15
\item $(\SO_n\times\GL_m,\omega_1\otimes\omega_1,\CC^n\otimes\CC^m)$
for $n\geq 3$, $\frac{1}{2}n\geq m\geq 1$.\\
We have $G_v^\circ\cong\SO_m\times\SO_{n-m}$ and $f(X)=\det(X^\top Q X)$
for $X\in\mat_{n,m}\cong\CC^n\otimes\CC^m$, $\deg(f)=2m$,
where $Q=g^\top Q g$ for $g\in\SO_n$.
% 16
\item $(\GL_1\times\Spin_7,\mu\otimes\mathrm{spinrep},\CC\otimes V^8)$.\\
We have $G_v^\circ\cong\G_2$ and $\deg(f)=2$, where $f$ is the relative 
invariant of SK I-15 for $m=1$, $n=8$.
% 17
\item $(\GL_2\times\Spin_7,\omega_1\otimes\mathrm{spinrep},\CC^2\otimes V^8)$.\\
We have $G_v^\circ\cong \SO_2\times\SL_3$ and $\deg(f)=4$, where $f$ is
the relative invariant of SK I-15 for $m=2$, $n=8$.
% 18
\item $(\GL_3\times\Spin_7,\omega_1\otimes\mathrm{spinrep},\CC^3\otimes V^8)$.\\
We have $G_v^\circ\cong \SO_3\times\SL_2$ and $\deg(f)=6$, where $f$ is
the relative invariant of SK I-15 for $m=3$, $n=8$.
% 19
\item $(\GL_1\times\Spin_9, \mu\otimes\mathrm{spinrep}, \CC\otimes V^{16})$.\\
We have $G_v^\circ\cong\Spin_7$ and
$\deg(f)=2$.
% 20
\item $(\GL_2\times\Spin_{10}, \omega_1\otimes\mathrm{halfspinrep},
\CC^2\otimes V^{16})$.\\
We have $G_v^\circ\cong \SL_2\times\G_2$ and $\deg(f)=4$.
% 21
\item $(\GL_3\times\Spin_{10}, \omega_1\otimes\mathrm{halfspinrep},
\CC^2\otimes V^{16})$.\\
We have $G_v^\circ\cong\SO_3\times\SL_2$ and $\deg(f)=12$.
% 22
\item $(\GL_1\times\Spin_{11},\mu\otimes\mathrm{spinrep},\CC\otimes V^{32})$.\\
We have $G_v^\circ\cong \SL_5$ and $\deg(f)=4$.
% 23
\item $(\GL_1\times\Spin_{12},\mu\otimes\mathrm{halfspinrep},\CC\otimes V^{32})$.\\
We have $G_v^\circ\cong \SL_6$ and $\deg(f)=4$.
% 24
\item $(\GL_1\times\Spin_{14},\mu\otimes\mathrm{halfspinrep},\CC\otimes V^{64})$.\\
We have $G_v^\circ\cong \G_2\times\G_2$ and $\deg(f)=8$.
% 25
\item $(\GL_1\times\G_2,\mu\otimes\omega_2,\CC\otimes V^7)$.\\
We have $G_v^\circ\cong\SL_3$ and $\deg(f)=2$, where $f$ is the relative
invariant of SK I-15 for $m=1$, $n=7$.
% 26
\item $(\GL_2\times\G_2,\omega_1\otimes\omega_2,\CC^2\otimes V^7)$.\\
We have $G_v^\circ\cong\GL_2$ and $\deg(f)=4$, where $f$ is the relative
invariant of SK I-15 for $m=2$, $n=7$.
% 27
\item $(\GL_1\times \E_6,\mu\otimes\omega_1, \CC\otimes V^{27})$.\\
We have $G_v^\circ\cong \F_4$ and $\deg(f)=4$.
% 28
\item $(\GL_2\times\E_6,\omega_1\otimes\omega_1,\CC^2\otimes V^{27})$.\\
We have $G_v^\circ\cong\SO_8$ and $\deg(f)=12$.
% 29
\item $(\GL_1\times\E_7,\mu\otimes\omega_6,\CC\otimes V^{56})$.\\
We have $G_v^\circ\cong \E_6$ and $\deg(f)=4$.
\end{enumerate}
% SK II
\item[{\rm \textbf{SK II}}] Non-regular irreducible reduced 
prehomogeneous modules with non-constant relative invariant.
\begin{enumerate}
\item $(\GL_1\times\Sp_n\times\SO_3,\mu\otimes\omega_1\otimes\omega_1,
\CC\otimes\CC^{2n}\otimes\CC^3)$.\\
We have $G_v^\circ\cong(\Sp_{n-2}\times\SO_2)\cdot \Un_{2n-3}$ and
$f(X)=\tr(X^\top J X Q)^2$ for $X\in\mat_{2n,3}\cong\CC\otimes\CC^{2n}\otimes\CC^3$.
\end{enumerate}
% SK III
\item[{\rm \textbf{SK III}}] Non-regular irreducible reduced 
prehomogeneous modules without non-constant relative invariants.
\begin{enumerate}
% 1 & 1'
\item $(G\times\GL_m,\rho\otimes\omega_1, V^n\otimes\CC^m)$,\\
where $\rho:G\Ra\GL(V^n)$ is an $n$-dimensional irreducible
representation of a semisimple algebraic group $G (\neq\SL_n)$
with $m>n\geq 3$.
We have $G_v^\circ\cong (G\times\GL_{m-n})\cdot\G^+_{n(m-n)}$.
The module $(G\times\SL_m,\rho\otimes\omega_1)$ is prehomogeneous
with $G_v^\circ\cong (G\times\SL_{m-n})\cdot\G^+_{n(m-n)}$.
% 2 % 2'
\item $(\SL_n\times\GL_m,\omega_1\otimes\omega_1, \CC^n\otimes\CC^m)$
for $\frac{1}{2}m\geq n\geq 1$.\\
We have $G_v^\circ\cong (\SL_n\times\GL_{m-n})\cdot\G^+_{n(m-n)}$.
The module $(\SL_n\times\SL_m,\omega_1\otimes\omega_1)$ is 
prehomogeneous with
$G_v^\circ\cong(\SL_n\times\SL_{m-n})\cdot\G^+_{n(m-n)}$.
% 3 & 3'
\item $(\GL_{2n+1},\omega_2,\bigwedge^2\CC^{2n+1})$ for $n\geq 2$.\\
We have $G_v^\circ\cong (\Sp_n\times\GL_1)\cdot\G^+_{2n}$.
The module $(\SL_{2n+1},\omega_2)$ is prehomogeneous with
$G_v^\circ\cong \Sp_n\cdot\G^+_{2n}$.
% 4 & 4'
\item $(\GL_2\times\SL_{2n+1},\omega_1\otimes\omega_2,
\CC^2\otimes\bigwedge^2\CC^{2n+1})$ for $n\geq 2$.\\
We have $G_v^\circ\cong(\GL_1\times\SL_2)\cdot\G^+_{2n}$
(see lemma 1.4 in Kimura et al. \cite{kimuraI}).
The module
$(\SL_2\times\SL_{2n+1},\omega_1\otimes\omega_2)$ is prehomogeneous with
$G_v^\circ\cong\SL_2\cdot\G^+_{2n}$.
% 5 & 5'
\item $(\Sp_n\times\GL_{2m+1},\omega_1\otimes\omega_1,
\CC^{2n}\otimes\CC^{2m+1})$ for $n>2m+1\geq 1$.\\
We have $G_v^\circ\cong (\GL_1\times\Sp_m\times\Sp_{n-m})\cdot\Un_{2n-1}$.
The module
$(\Sp_n\times\SL_{2m+1},\omega_1\otimes\omega_1)$ is prehomogeneous
with $G_v^\circ\cong (\Sp_m\times\Sp_{n-m})\cdot\Un_{2n-1}$.
% 6 & 6'
\item $(\GL_1\times\Spin_{10},\mu\otimes\mathrm{halfspinrep},\CC\otimes V^{16})$.\\
We have $G_v^\circ\cong(\GL_1\times\Spin_7)\cdot\G^+_8$.
The module
$(\Spin_{10},\mathrm{halfspinrep})$ is prehomogeneous with
$G_v^\circ\cong\Spin_7\cdot\G^+_8$.
\end{enumerate}
\end{itemize}

\subsection{Non-irreducible simple prehomogeneous modules}

The simple prehomogeneous, including the non-irreducible ones,
were classified by Kimura, thus we will label them
by Ks $n$, where $n$ is the number of the module in 
\S 3 of Kimura's article \cite{kimuraS}.
\index{Ks}

In this section, it is understood that each representation $\rho_i$
of the simple group is composed with a scalar multiplication $\mu$ of
$\GL_1^k$. We shall simply write $\rho_i$ instead of
$\mu\otimes\rho_i$. In some cases, a module
$V_1\oplus\ldots\oplus V_k$ will be prehomogeneous even with fewer
than $k$ scalar multiplications, in which case we will state this
fact explicitely. We shall also state the connected component $G_v^\circ$
of the generic isotropy subgroup and the relative invariants
$f_1,\ldots,f_l$ where they exist.
\index{prehomogeneous module!simple}
\index{prehomogeneous module!classification}

Let $G=\GL_1^k\times G_{\mathrm{s}}$ be a reductive group, where
$G_{\mathrm{s}}$ is a simple algebraic group, let
$(\rho_1,V_1)$, \ldots, $(\rho_k,V_k)$ be irreducible
$G_{\mathrm{s}}$-modules
and $V=V_1\oplus\ldots\oplus V_k$ a $G_{\mathrm{s}}$-module with representation $\rho=\rho_1\oplus\ldots\oplus\rho_k$.
Then $(G,\rho,V)$ is  equivalent to one of the following:
\begin{itemize}
% Ks I
\item[{\rm \textbf{Ks I}}] Regular non-irreducible
simple prehomogeneous modules.
\begin{enumerate}
% 1
\item $(\GL_1^2\times\SL_n,\omega_1\oplus\omega_1^*,\CC^n\oplus\CC^{n*})$ for $n\geq3$.\\
We have $G_v^\circ\cong \GL_1\times\SL_{n-1}$ and
$f_1(x,y)=\lag x|y\rag$, where $(x,y)\in\CC^n\oplus\CC^{n*}$
and $\lag\cdot\mid\cdot\rag$ is the dual pairing.
The module $(\GL_1\times\SL_n,(\mu\otimes\omega_1)\oplus\omega_1^*)$ is prehomogeneous with $G_v^\circ=\SL_{n-1}$.
% 2
\item $(\GL_1^n\times\SL_n,
\omega_1^{\oplus n},(\CC^n)^{\oplus n})$ for $n\geq2$.\\
We have $G_v^\circ\cong \GL_1^{n-1}$ and
$f_1(X)=\det(X)$ for $X\in\mat_n\cong(\CC^n)^{\oplus n}$.
The module $(\GL_1\times\SL_n,\mu\otimes\omega_1^{\oplus n})$ is prehomogeneous with $G_v^\circ=\{1\}$.
% 3
\item $(\GL_1^{n+1}\times\SL_n,
\omega_1^{\oplus n+1},(\CC^n)^{\oplus n+1})$ for $n\geq2$.\\
We have $G_v^\circ\cong \{1\}$ and
$f_i(X)=\det(x_1,\ldots,\cancel{x_i},\ldots,x_{n+1})$ for
$X=(x_1,\ldots,x_{n+1})\in\mat_{n,n+1}\cong(\CC^n)^{\oplus n+1}$.
% 4
\item $(\GL_1^{n+1}\times\SL_n,
\omega_1^{\oplus n}\oplus\omega_1^*,(\CC^n)^{\oplus n}\oplus\CC^{n*})$ for $n\geq3$.\\
We have $G_v^\circ\cong \{1\}$ and $f_1(X)=\lag x_1|y\rag, \ldots,
f_n(X)=\lag x_n|y\rag$,
$f_{n+1}(X)=\det(x_1,\ldots,x_n)$ for
$X=(x_1,\ldots,x_n,y)\in(\CC^n)^{\oplus n}\oplus\CC^{n*}$.
% 5
\item $(\GL_1^3\times\SL_{2n},\omega_2\oplus\omega_1\oplus\omega_1,
\bigwedge^2 \CC^{2n}\oplus\CC^{2n}\oplus\CC^{2n})$ for $n\geq2$.\\
We have $G_v^\circ\cong\GL_1\times\Sp_{n-1}$ and $f_1(X,y,z)=\Pf(X)$,
$f_2(X,y,z)=y^\top X^\# z$, where
$(X,y,z)\in\bigwedge^2 \CC^{2n}\oplus\CC^{2n}\oplus\CC^{2n}$ and
$X^\#$ is the cofactor matrix of $X$.
The module $(\GL^2_1\times\SL_{2n},(\mu\otimes\omega_2)\oplus(\mu\otimes(\omega_1\oplus\omega_1)))$ is prehomogeneous with
$G_v^\circ\cong\Sp_{n-1}$.
% 6
\item $(\GL_1^3\times\SL_{2n},\omega_2\oplus\omega_1\oplus\omega_1^*,
\bigwedge^2 \CC^{2n}\oplus\CC^{2n}\oplus\CC^{2n*})$ for $n\geq2$.\\
We have $G_v^\circ\cong\GL_1\times\Sp_{n-1}$ and $f_1(X,y,z)=\Pf(X)$,
$f_2(X,y,z)=\lag y|z\rag$, where
$(X,y,z)\in\bigwedge^2 \CC^{2n}\oplus\CC^{2n}\oplus\CC^{2n*}$.
The module
$(\GL^2_1\times\SL_{2n},(\mu\otimes\omega_2)\oplus(\mu\otimes(\omega_1\oplus\omega_1^*))$
is prehomogeneous with $G_v^\circ\cong\Sp_{n-1}$.
% 7
\item $(\GL_1^3\times\SL_{2n},\omega_2\oplus\omega_1^*\oplus\omega_1^*,
\bigwedge^2 \CC^{2n}\oplus\CC^{2n*}\oplus\CC^{2n*})$ for $n\geq3$.\\
We have $G_v^\circ\cong\GL_1\times\Sp_{n-1}$ and $f_1(X,y,z)=\Pf(X)$,
$f_2(X,y,z)=y^\top X z$, where
$(X,y,z)\in\bigwedge^2 \CC^{2n}\oplus\CC^{2n*}\oplus\CC^{2n*}$.
The module $(\GL^2_1\times\SL_{2n},(\mu\otimes\omega_2)\oplus(\mu\otimes(\omega_1^*\oplus\omega_1^*))$ is prehomogeneous with
$G_v^\circ\cong\Sp_{n-1}$.
% 8
\item $(\GL_1^2\times\SL_{2n+1},\omega_2\oplus\omega_1,
\bigwedge^2 \CC^{2n+1}\oplus\CC^{2n+1})$ for $n\geq2$.\\
We have $G_v^\circ\cong \GL_1\times\Sp_{n}$.
The module
$(\GL_1\times\Sp_{2n+1},\mu\otimes(\omega_2\oplus\omega_1))$
is prehomogeneous with $G_v^\circ\cong\Sp_n$, see p. 94 in \cite{kimuraS} for the relative invariant.
% 9
\item $(\GL_1^4\times\SL_{2n+1},\omega_2\oplus\omega_1\oplus\omega_1\oplus\omega_1,\bigwedge^2\CC^{2n+1}\oplus\CC^{2n+1}\oplus\CC^{2n+1}\oplus\CC^{2n+1})$ for $n\geq2$.\\
We have $G_v^\circ\cong\Sp_{n-1}$, see p. 94 in \cite{kimuraS} for the relative invariants.
% 10
\item $(\GL_1^4\times\SL_{2n+1},\omega_2\oplus\omega_1\oplus\omega_1^*\oplus\omega_1^*,\bigwedge^2\CC^{2n+1}\oplus\CC^{2n+1}\oplus\CC^{2n+1*}\oplus\CC^{2n+1*})$ for $n\geq2$.\\
We have $G_v^\circ\cong \Sp_{n-1}$ and
$f_1(X)=\Pf\begin{pmatrix} A & x \\ x^\top & 0 \end{pmatrix}$,
$f_2(X)=\lag x|y\rag$, $f_3(X)=\lag x|z\rag$, $f_4= y^\top A z$
for $X=(A,x,y,z)\in\bigwedge^2\CC^{2n+1}\oplus\CC^{2n+1}\oplus\CC^{2n+1*}\oplus\CC^{2n+1*}$.
% 11
\item $(\GL_1^2\times\SL_n,2\omega_1\oplus\omega_1,\Sym^2\CC^n\oplus\CC^n)$ for $n\geq2$.\\
We have $G_v^\circ\cong \SO_{n-1}$ and $f_1(X)=\det(A)$,
$f_2(X)=x^\top A^\# x$ for $X=(A,x)\in\Sym^2\CC^n\oplus\CC^n$.
% 12
\item $(\GL_1^2\times\SL_n,2\omega_1\oplus\omega_1^*,\Sym^2\CC^n\oplus\CC^{n*})$ for $n\geq3$.\\
We have $G_v^\circ\cong \SO_{n-1}$ and $f_1(X)=\det(A)$,
$f_2(X)=x^\top A x$ for $X=(A,x)\in\Sym^2\CC^n\oplus\CC^n$.
% 13
\item $(\GL_1^2\times\SL_7,\omega_3\oplus\omega_1,\bigwedge^3\CC^7\oplus\CC^7)$.\\
We have $G_v^\circ\cong\SL_3$, see p. 96 in \cite{kimuraS} for the relative
invariants.
% 14
\item $(\GL_1^2\times\SL_7,\omega_3\oplus\omega_1^*,\bigwedge^3\CC^7\oplus\CC^{7*})$.\\
We have $G_v^\circ\cong\SL_3$, see p. 96 in \cite{kimuraS} for the relative
invariants.
% 15
\item $(\GL_1^2\times\Spin_8,\mathrm{spinrep}\oplus\mathrm{halfspinrep},V^8\oplus V^8)$.\\
We have $G_v^\circ\cong\G_2$ and two quadratic invariants
$f_1(x,y)=q_1(x)$, $f_2(x,y)=q_2(y)$ for $(x,y)\in V^8\oplus V^8$.
% 16
\item $(\GL_1^2\times\Spin_7,
\mathrm{vecrep}\oplus\mathrm{spinrep},V^7\oplus V^8)$.\\
We have $G_v^\circ\cong\SL_3$ and two quadratic invariants
$f_1(x,y)=q_1(x)$, $f_2(x,y)=q_2(y)$ for $(x,y)\in V^7\oplus V^8$.
% 17
\item $(\GL_1^2\times\Spin_{10},\mathrm{halfspinrep}_\mathrm{even}\oplus\mathrm{halfspinrep}_\mathrm{even},
V^{16}\oplus V^{16})$.\\
We have $G_v^\circ\cong\GL_1\times\G_2$, see p. 96 in \cite{kimuraS}
for the relative invariants. The module
$(\GL_1\times\Spin_{10},\mu\otimes(\mathrm{halfspinrep}_\mathrm{even}\oplus\mathrm{halfspinrep}_\mathrm{even}))$
is prehomogeneous with $G_v^\circ\cong\G_2$.
% 18
\item $(\GL_1^2\times\Spin_{10},\mathrm{vecrep}\oplus\mathrm{halfspinrep}, V^{10}\oplus V^{16})$.\\
We have $G_v^\circ\cong\Spin_7$, see p. 97 in \cite{kimuraS}
for the relative invariants.
% 19
\item $(\GL_1^2\times\Spin_{12},\mathrm{vecrep}\oplus\mathrm{halfspinrep}, V^{12}\oplus V^{32})$.\\
We have $G_v^\circ\cong\SL_5$, see p. 97 in \cite{kimuraS}
for the relative invariants.
% 20
\item $(\GL_1^2\times\Sp_{n},\omega_1\oplus\omega_1,
\CC^{2n}\oplus \CC^{2n})$ for $n\geq2$.\\
We have $G_v^\circ\cong\GL_1\times\Sp_{n-1}$, see p. 97 in \cite{kimuraS}
for the relative invariant. The module
$(\GL_1\times\Sp_{n},\mu\otimes(\omega_1\oplus\omega_1))$ is
prehomogeneous with $G_v^\circ\cong\Sp_{n-1}$.
% 21
\item $(\GL_1^2\times\Sp_{3},\omega_3\oplus\omega_1,
V^{14}\oplus \CC^6)$.\\
We have $G_v^\circ\cong\SL_2$, see p. 97 in \cite{kimuraS}
for the relative invariants.
\end{enumerate}
% Ks II
\item[{\rm \textbf{Ks II}}] Non-regular non-irreducible simple
prehomogeneous modules.
\begin{enumerate}
% 1
\item $(\GL_1^k\times\SL_n,\omega_1^{\oplus k},(\CC^n)^{\oplus k})$
for $2\leq k\leq n-1$.\\
We have $G_v^\circ\cong (\GL_1^k\times\SL_{n-k})\cdot \G^+_{k(n-k)}$.
The module $(\SL_n,\omega_1^{\oplus k})$ is
prehomogeneous with
$G_v^\circ\cong \SL_{n-k}\cdot \G^+_{k(n-k)}$.
% 2
\item $(\GL_1^k\times\SL_n,\omega_1^{\oplus k-1}\oplus\omega_1^*,
(\CC^n)^{\oplus k-1}\oplus\CC^{n*})$ for $3\leq k\leq n$.\\
We have $G_v^\circ\cong (\GL_1\times\SL_{n-k+1})\cdot \G^+_{(n-k+1)(k-2)}$ and $f_1(X)=\lag x_1|y\rag$, \ldots,
$f_{k-1}(X)=\lag x_{k-1}|y\rag$ for $X=(x_1,\ldots,x_{k-1},y)\in
(\CC^n)^{\oplus k-1}\oplus\CC^{n*}$.
The module $(\GL_1^{k-1}\times\SL_n,(\mu\otimes\omega_1^{\oplus k-1})\oplus\omega_1^*)$
is prehomogeneous with
$G_v^\circ\cong \SL_{n-k+1}\cdot \G^+_{(n-k+1)(k-2)}$.
% 3
\item $(\GL_1^2\times\SL_{2n+1},\omega_2\oplus\omega_2,
\bigwedge^2\CC^{2n+1}\oplus\bigwedge^2\CC^{2n+1})$ for $n\geq 2$.\\
We have $G_v^\circ\cong \GL_1^2\cdot\G_{2n}^+$. The module
$(\SL_{2n+1},\omega_2\oplus\omega_2)$ is prehomogeneous with
$G_v^\circ\cong \G^+_{2n}$.
% 4
\item $(\GL_1^2\times\SL_{2n},\omega_2\oplus\omega_1,
\bigwedge^2\CC^{2n}\oplus\CC^{2n})$ for $n\geq 2$.\\
We have $G_v^\circ\cong (\GL_1\times\Sp_{n-1})\cdot\Un_{n-1}$ and
$f_1(X)=\Pf(X)$ where $X\in\bigwedge^2\CC^{2n}$. The module
$(\GL_1\times\SL_{2n},\mu\otimes(\omega_2\oplus\omega_1))$
is prehomogeneous with $G_v^\circ\cong \Sp_{n-1}\cdot\Un_{n-1}$.
% 5
\item $(\GL_1^2\times\SL_{2n},\omega_2\oplus\omega_1^*,
\bigwedge^2\CC^{2n}\oplus\CC^{2n*})$ for $n\geq 3$.\\
We have $G_v^\circ\cong (\GL_1\times\Sp_{n-1})\cdot\Un_{n-1}$ and
$f_1(X)=\Pf(X)$ where $X\in\bigwedge^2\CC^{2n}$. The module
$(\GL_1\times\SL_{2n},\mu\otimes(\omega_2\oplus\omega_1^*))$
is prehomogeneous with $G_v^\circ\cong \Sp_{n-1}\cdot\Un_{n-1}$.
% 6
\item $(\GL_1^4\times\SL_{2n},\omega_2\oplus\omega_1\oplus\omega_1\oplus\omega_1,\bigwedge^2\CC^{2n}\oplus\CC^{2n}\oplus\CC^{2n}\oplus\CC^{2n})$ for $n\geq 2$.\\
We have $G_v^\circ\cong \Sp_{n-2}\cdot\Un_{2n-3}$ and $f_1(X)=\Pf(A)$,
$f_2(X)=x^\top A^\# y$, $f_3(X)=y^\top A^\# z$, $f_4(X)=z^\top A x$
for $X=(A,x,y,z)\in\bigwedge^2\CC^{2n}\oplus\CC^{2n}\oplus\CC^{2n}\oplus\CC^{2n}$.
% 7
\item $(\GL_1^4\times\SL_{2n},\omega_2\oplus\omega_1\oplus\omega_1\oplus\omega_1^*,\bigwedge^2\CC^{2n}\oplus\CC^{2n}\oplus\CC^{2n}\oplus\CC^{2n*})$ for $n\geq 2$.\\
We have $G_v^\circ\cong \Sp_{n-2}\cdot\Un_{2n-3}$ and $f_1(X)=\Pf(A)$,
$f_2(X)=x^\top A^\# y$, $f_3(X)=\lag x|z\rag$, $f_4(X)=\lag y|z\rag$
for $X=(A,x,y,z)\in\bigwedge^2\CC^{2n}\oplus\CC^{2n}\oplus\CC^{2n}\oplus\CC^{2n*}$.
% 8
\item $(\GL_1^4\times\SL_{2n},\omega_2\oplus\omega_1\oplus\omega_1^*\oplus\omega_1^*,\bigwedge^2\CC^{2n}\oplus\CC^{2n}\oplus\CC^{2n*}\oplus\CC^{2n*})$ for $n\geq 3$.\\
We have $G_v^\circ\cong \Sp_{n-2}\cdot\Un_{2n-3}$ and $f_1(X)=\Pf(A)$,
$f_2(X)=\lag x|y\rag$, $f_3(X)=\lag x|z\rag$, $f_4(X)=y^\top A z$
for $X=(A,x,y,z)\in\bigwedge^2\CC^{2n}\oplus\CC^{2n}\oplus\CC^{2n*}\oplus\CC^{2n*}$.
% 9
\item $(\GL_1^4\times\SL_{2n},\omega_2\oplus\omega_1^*\oplus\omega_1^*\oplus\omega_1^*,\bigwedge^2\CC^{2n}\oplus\CC^{2n*}\oplus\CC^{2n*}\oplus\CC^{2n*})$ for $n\geq 3$.\\
We have $G_v^\circ\cong \Sp_{n-2}\cdot\Un_{2n-3}$ and $f_1(X)=\Pf(A)$,
$f_2(X)=x^\top A y$, $f_3(X)=y^\top A z$, $f_4(X)=z^\top A x$
for $X=(A,x,y,z)\in\bigwedge^2\CC^{2n}\oplus\CC^{2n*}\oplus\CC^{2n*}\oplus\CC^{2n*}$.
% 10
\item $(\GL_1^2\times\SL_{2n+1},\omega_2\oplus\omega_1^*,
\bigwedge^2\CC^{2n+1}\oplus\CC^{2n+1*})$ for $n\geq 2$.\\
We have $G_v^\circ\cong(\GL_1^2\times\Sp_{n-1})\cdot\Un_{4n-2}$. The
module $(\SL_{2n+1},\omega_2\oplus\omega_1^*)$ is prehomogeneous with
$G_v^\circ\cong\Sp_{n-1}\cdot\Un_{4n-2}$.
% 11
\item $(\GL_1^3\times\SL_{2n+1},\omega_2\oplus\omega_1\oplus\omega_1,
\bigwedge^2\CC^{2n+1}\oplus\CC^{2n+1}\oplus\CC^{2n+1})$ for $n\geq 2$.\\
We have $G_v^\circ\cong (\GL_1\times\Sp_{n-1})\cdot\Un_{2n-1}$, see
p. 99 in \cite{kimuraS} for the relative invariants.
The module $(\GL_1^2\times\SL_{2n+1},(\mu\otimes(\omega_2\oplus\omega_1))\oplus(\mu\otimes\omega_1))$ is prehomogeneous with
$G_v^\circ\cong\Sp_{n-1}\cdot\Un_{2n-1}$.
% 12
\item $(\GL_1^3\times\SL_{2n+1},\omega_2\oplus\omega_1\oplus\omega_1^*,
\bigwedge^2\CC^{2n+1}\oplus\CC^{2n+1}\oplus\CC^{2n+1*})$ for $n\geq 2$.\\
We have $G_v^\circ\cong (\GL_1\times\Sp_{n-1})\cdot\Un_{2n-1}$, see
p. 99 in \cite{kimuraS} for the relative invariants.
The module $(\GL_1^2\times\SL_{2n+1},(\mu\otimes(\omega_2\oplus\omega_1))\oplus(\mu\otimes\omega_1^*))$ is prehomogeneous with
$G_v^\circ\cong\Sp_{n-1}\cdot\Un_{2n-1}$.
% 13
\item $(\GL_1^3\times\SL_{2n+1},\omega_2\oplus\omega_1^*\oplus\omega_1^*,
\bigwedge^2\CC^{2n+1}\oplus\CC^{2n+1*}\oplus\CC^{2n+1*})$ for $n\geq 2$.\\
We have $G_v^\circ\cong (\GL_1\times\Sp_{n-1})\cdot\Un_{2n-1}$ and
$f_1(X)=x^\top A y$ for $X=(A,x,y)\in\bigwedge^2\CC^{2n+1}\oplus\CC^{2n+1*}\oplus\CC^{2n+1*}$.
The module $(\GL_1\times\SL_{2n+1},\omega_2\oplus(\mu\otimes(\omega_1^*\oplus\omega_1^*))$ is prehomogeneous with
$G_v^\circ\cong\Sp_{n-1}\cdot\Un_{2n-2}$.
% 14
\item $(\GL_1^4\times\SL_{2n+1},\omega_2\oplus\omega_1^*\oplus\omega_1^*\oplus\omega_1^*,\bigwedge^2\CC^{2n+1}\oplus\CC^{2n+1}\oplus\CC^{2n+1}\oplus\CC^{2n+1})$ for $n\geq 2$.\\
We have $G_v^\circ\cong(\GL_1\times\Sp_{n-2})\cdot\Un_{4n-6}$ and
$f_1(X)=x^\top A y$, $f_2(X)=y^\top A z$, $f_3(X)=z^\top A x$.
The module $(\GL_1^3\times\SL_{2n+1},\omega_2\oplus(\mu\otimes\omega_1^*)\oplus(\mu\otimes\omega_1^*)\oplus(\mu\otimes\omega_1^*))$ is prehomogeneous with $G_v^\circ\cong\Sp_{n-2}\cdot\Un_{4n-6}$.
% 15
\item $(\GL_1^2\times\SL_6,\omega_3\oplus\omega_1,\bigwedge^3\CC^6\oplus\CC^6)$.\\
We have $G_v^\circ\cong(\GL_1\times\SL_2\times\SL_2)\cdot\G^+_4$,
see p. 100 in \cite{kimuraS} for the relative invariant.
The module $(\GL_1\times\SL_6,\mu\otimes(\omega_3\oplus\omega_1))$ is
prehomogeneous with $G_v^\circ\cong(\SL_2\times\SL_2)\cdot\G^+_4$.
% 16
\item $(\GL_1^3\times\SL_6,\omega_3\oplus\omega_1\oplus\omega_1,\bigwedge^3\CC^6\oplus\CC^6\oplus\CC^6)$.\\
We have $G_v^\circ\cong\GL_1^2\cdot\G^+_4$,
see p. 100 in \cite{kimuraS} for the relative invariant.
The module $(\GL_1\times\SL_6,\mu\otimes(\omega_3\oplus\omega_1\oplus\omega_1))$ is
prehomogeneous with $G_v^\circ\cong\G^+_4$.
% 17
\item $(\GL_1^3\times\Sp_n,\omega_1\oplus\omega_1\oplus\omega_1,
\CC^{2n}\oplus\CC^{2n}\oplus\CC^{2n})$ for $n\geq 2$.\\
We have $G_v^\circ\cong \Sp_{n-2}\cdot\Un_{2n-3}$, see p. 100 in \cite{kimuraS} for the relative invariants.
% 18
\item $(\GL_1^2\times\Sp_2,\omega_2\oplus\omega_1,
V^5\oplus\CC^4)$.\\
We have $G_v^\circ\cong \GL_1\cdot\Un_2$, see p. 100 in \cite{kimuraS} for the relative invariant. The module $(\GL_1\times\Sp_2,(\mu\otimes\omega_2)\oplus\omega_1)$ is prehomogeneous with $G_v^\circ\cong\Un_2$.
% 19
\item $(\GL_1^3\times\SL_5,\omega_2\oplus\omega_2\oplus\omega_1^*,
\bigwedge^2\CC^5\oplus\bigwedge^2\CC^5\oplus\CC^5)$.\\
See proposition 1.1 in \cite{kimuraI}.
\end{enumerate}
\end{itemize}

\subsection{Two-simple prehomogeneous modules of type I}

In this and the following chapter we shall give a classification
of the non-irreducible two-simple prehomogeneous modules,
i.e. modules of the form
\begin{align*}
\bigl(&\GL_1^l\times G_1 \times G_2,\\
&(\rho_1 \otimes \tilde{\rho}_1)\oplus\ldots\oplus(\rho_k\otimes\tilde{\rho}_k)\
\oplus\
(\sigma_1\otimes 1)\oplus\ldots\oplus(\sigma_s\otimes 1)\
\oplus\
(1\otimes\tau_1)\oplus\ldots\oplus(1\otimes\tau_t),\\
&V_1\oplus\ldots \oplus V_l \bigr),
\end{align*}
where $G_1$ and $G_2$ are simple algebraic groups, $l=k+s+t$, and
the $\rho_i$, $\sigma_j$ (resp.~$\tilde{\rho}_i$, $\tau_j$) are
irreducible representations of $G_1$ (resp.~$G_2$).
As in the previous chapter, it is understood that each of these
representations is composed with a scalar multiplication of $\GL_1^k$.
First, we give the classification of the type I-modules,
i.e. at least one of the modules
$(\GL_1\times G_1\times G_2, \rho_i\otimes\tilde{\rho}_i)$
is a non-trivial prehomogeneous module.
These were classified by Kimura et al. \cite{kimuraI}, thus we
shall refer to them as KI $n$, where $n$ is the number of the module
in \S 3 of \cite{kimuraI}. We shall state the non-irreducible modules
only, as the irreducible ones already appear in the table
SK or as castling transformas of those (see also
theorem 1.5 in \cite{kimuraI}).
In the next chapter, we shall classify the remaining two-simple modules
of type II.

Let $(G,\rho,V)$ be a two-simple prehomogeneous module of type I.
Then it is  equivalent to one of the following:
\begin{itemize}
% KI I
\item[{\rm \textbf{KI I}}] Regular two-simple prehomogeneous modules
of type I.
\begin{enumerate}
% 1
\item $(\GL_1^2\times\SL_4\times\SL_2,
(\omega_2\otimes\omega_1)\oplus(\omega_1\otimes\omega_1))$.\\
We have $G_v^\circ\cong\{1\}$.
% 2
\item $(\GL_1^3\times\SL_4\times\SL_2,
(\omega_2\otimes\omega_1)\oplus(\omega_1\otimes 1)\oplus(\omega_1\otimes 1))$.\\
We have $G_v^\circ\cong \GL_1$.
The module $(\GL_1^2\times\SL_4\times\SL_2,
(\omega_2\otimes\omega_1)\oplus(\omega_1\otimes 1)\oplus(\omega_1\otimes 1))$ is prehomogeneous with $G_v^\circ\cong\{1\}$.
% 3
\item $(\GL_1^2\times\SL_4\times\SL_3,
(\omega_2\otimes\omega_1)\oplus(\omega_1\otimes 1))$.\\
We have $G_v^\circ\cong \SO_3$.
% 4
\item $(\GL_1^3\times\SL_4\times\SL_3,
(\omega_2\otimes\omega_1)\oplus(\omega_1\otimes 1)\oplus(1\otimes\omega_1^{(*)}))$.\\
We have $G_v^\circ\cong\SO_2$.
% 5
\item $(\GL_1^3\times\SL_4\times\SL_4,
(\omega_2\otimes\omega_1)\oplus(\omega_1\otimes 1)\oplus(1\otimes\omega_1^*))$.\\
We have $G_v^\circ\cong\SO_2$.
% 6
\item $(\GL_1^3\times\SL_5\times\SL_2,
(\omega_2\otimes\omega_1)\oplus(\omega_1^*\otimes 1)\oplus(\omega_1^{(*)}\otimes1))$.\\
We have $G_v^\circ\cong\{1\}$.
% 7
\item $(\GL_1^2\times\SL_5\times\SL_3,
(\omega_2\otimes\omega_1)\oplus(1\otimes\omega_1^{(*)}))$.\\
We have $G_v^\circ\cong\SO_2$.
% 8
\item $(\GL_1^2\times\SL_5\times\SL_8,
(\omega_2\otimes\omega_1)\oplus(1\otimes\omega_1^*))$.\\
We have $G_v^\circ\cong\SO_2$.
% 9
\item $(\GL_1^2\times\SL_5\times\SL_9,
(\omega_2\otimes\omega_1)\oplus(1\otimes\omega_1^*))$.\\
We have $G_v^\circ\cong\GL_1\times\SL_2\times\SL_2$.
The module $(\GL_1\times\SL_5\times\SL_9,
(\omega_2\otimes\omega_1)\oplus(1\otimes\omega_1^*))$ is
prehomogeneous with $G_v^\circ\cong\SL_2\times\SL_2$.
% 10
\item $(\GL_1^3\times\Sp_n\times\SL_{2m},
(\omega_1\otimes\omega_1)\oplus(1\otimes\omega_1^{(*)})\oplus(1\otimes\omega_1^{(*)}))$.\\
We have $G_v^\circ\cong\GL_1\times\Sp_{n-m}\times\Sp_{m-1}$. The module
$(\GL_1^2\times\Sp_n\times\Sp_{2m},
(\omega_1\otimes\omega_1)\oplus(1\otimes\omega_1^{(*)})\oplus(1\otimes\omega_1^{(*)}))$ is prehomogeneous with $G_v^\circ\cong \Sp_{n-m}\times\Sp_{m-1}$.
% 11
\item $(\GL_1^2\times\Sp_n\times\SL_2,(\omega_1\otimes\omega_1)\oplus(1\otimes 2\omega_1))$.\\
We have $G_v^\circ\cong\Sp_{n-1}\times\SO_2$.
% 12
\item $(\GL_1^2\times\Sp_n\times\SL_3,(\omega_1\otimes\omega_1)\oplus(1\otimes 3\omega_1))$.\\
We have $G_v^\circ\cong \Sp_{n-1}$.
% 13
\item $(\GL_1^3\times\Sp_n\times\SL_2,(\omega_1\otimes\omega_1)\oplus(1\otimes2\omega_2)\oplus(1\otimes\omega_1))$.
We have $G_v^\circ\cong \Sp_{n-1}$.
% 14
\item $(\GL_1^2\times\Sp_n\times\SL_{2m+1},(\omega_1\otimes\omega_1)\oplus(\omega_1\otimes 1))$.\\
We have $G_v^\circ\cong \GL_1\times\Sp_m\times\Sp_{n-m-1}$. The module
$(\GL_1\times\Sp_n\times\SL_{2m+1},(\omega_1\otimes\omega_1)\oplus(\omega_1\otimes 1))$ is prehomogeneous with $G_v^\circ\cong\Sp_m\times\Sp_{n-m-1}$.
% 15
\item $(\GL_1^4\times\Sp_n\times\SL_{2m+1},(\omega_1\otimes\omega_1)\oplus(\omega_1\otimes 1)\oplus(1\otimes(\omega_1\oplus\omega_1)^{(*)}))$.\\
We have $G_v^\circ\cong\Sp_{m-1}\times\Sp_{n-m-1}$.
% 16
\item $(\GL_1^3\times\Sp_2\times\SL_3,(\omega_1\otimes\omega_1)\oplus(\omega_2\otimes 1)\oplus(1\otimes\omega_1^*))$.\\
We have $G_v^\circ\cong\GL_1$.
The module $(\GL_1^2\times\Sp_2\times\SL_3,(\omega_1\otimes\omega_1)\oplus(\omega_2\otimes 1)\oplus(1\otimes\omega_1^*))$ is
prehomogeneous with $G_v^\circ\cong\{1\}$.
% 17
\item $(\GL_1^2\times\Sp_2\times\SL_2,(\omega_2\otimes\omega_1)\oplus(\omega_1\otimes 1))$.\\
We have $G_v^\circ\cong\SO_2$.
% 18
\item $(\GL_1^3\times\Sp_2\times\SL_2,(\omega_2\otimes\omega_1)\oplus(\omega_1\otimes 1)\oplus(1\otimes\omega_1))$.\\
We have $G_v^\circ\cong\{1\}$.
% 19
\item $(\GL_1^3\times\Sp_2\times\SL_4,(\omega_2\otimes\omega_1)\oplus(\omega_1\otimes 1)\oplus(1\otimes\omega_1^*))$.\\
We have $G_v^\circ\cong\{1\}$.
% 20
\item $(\GL_1^2\times\SO_n\times\SL_m,(\omega_1\otimes\omega_1)\oplus(1\otimes\omega_1^{(*)}))$.\\
We have $G_v^\circ\cong \SO_{m-1}\times\SO_{n-m}$.
% 21
\item $(\GL_1^2\times\Spin_7\times\SL_2,
(\mathrm{spinrep}\otimes\omega_1)\oplus(1\otimes\omega_1))$.\\
We have $G_v^\circ\cong\SL_3$.
% 22
\item $(\GL_1^2\times\Spin_7\times\SL_3,
(\mathrm{spinrep}\otimes\omega_1)\oplus(1\otimes\omega_1^{(*)}))$.\\
We have $G_v^\circ\cong\SL_2\times\SO_2$.
% 23
\item $(\GL_1^2\times\Spin_7\times\SL_6,
(\mathrm{spinrep}\otimes\omega_1)\oplus(1\otimes\omega_1^*))$.\\
We have $G_v^\circ\cong\SL_2\times\SO_2$.
% 24
\item $(\GL_1^2\times\Spin_7\times\SL_7,
(\mathrm{spinrep}\otimes\omega_1)\oplus(1\otimes\omega_1^*))$.\\
We have $G_v^\circ\cong\SL_3$.
% 25
\item $(\GL_1^2\times\Spin_7\times\SL_2,
(\mathrm{vecrep}\otimes\omega_1)\oplus(\mathrm{spinrep}\otimes 1))$.\\
We have $G_v^\circ\cong\GL_2$.
% 26
\item $(\GL_1^3\times\Spin_7\times\SL_2,
(\mathrm{vecrep}\otimes\omega_1)\oplus(\mathrm{spinrep}\otimes 1)\oplus(1\otimes\omega_1))$.\\
We have $G_v^\circ\cong\SL_2$.
% 27
\item $(\GL_1^3\times\Spin_7\times\SL_6,
(\mathrm{vecrep}\otimes\omega_1)\oplus(\mathrm{spinrep}\otimes 1)\oplus(1\otimes\omega_1^*))$.\\
We have $G_v^\circ\cong\SL_2$.
% 28
\item $(\GL_1^2\times\Spin_8\times\SL_2,(\mathrm{vecrep}\otimes\omega_1)\oplus(\mathrm{halfspinrep}\otimes 1))$.\\
We have $G_v^\circ\cong\SL_3\times\SO_2$.
% 29
\item $(\GL_1^2\times\Spin_8\times\SL_3,(\mathrm{vecrep}\otimes\omega_1)\oplus(\mathrm{halfspinrep}\otimes 1))$.\\
We have $G_v^\circ\cong\SL_2\times\SO_3$.
% 30
\item $(\GL_1^3\times\Spin_8\times\SL_2,(\mathrm{vecrep}\otimes\omega_1)\oplus(\mathrm{halfspinrep}\otimes 1)\oplus(1\otimes\omega_1))$.\\
We have $G_v^\circ\cong\SL_3$.
% 31
\item $(\GL_1^3\times\Spin_8\times\SL_3,(\mathrm{vecrep}\otimes\omega_1)\oplus(\mathrm{halfspinrep}\otimes 1)\oplus(1\otimes\omega_1^{(*)}))$.\\
We have $G_v^\circ\cong\SL_2\times\SO_2$.
% 32
\item $(\GL_1^3\times\Spin_8\times\SL_6,(\mathrm{vecrep}\otimes\omega_1)\oplus(\mathrm{halfspinrep}\otimes 1)\oplus(1\otimes\omega_1^*))$.\\
We have $G_v^\circ\cong\SL_2\times\SO_2$.
% 33
\item $(\GL_1^3\times\Spin_8\times\SL_7,(\mathrm{vecrep}\otimes\omega_1)\oplus(\mathrm{halfspinrep}\otimes 1)\oplus(1\otimes\omega_1^*))$.\\
We have $G_v^\circ\cong\SL_3$.
% 34
\item $(\GL_1^2\times\Spin_{10}\times\SL_2,(\mathrm{halfspinrep}\otimes\omega_1)\oplus(1\otimes2\omega_1))$.\\
We have $G_v^\circ\cong\G_2\times\SO_3$.
% 35
\item $(\GL_1^2\times\Spin_{10}\times\SL_2,(\mathrm{halfspinrep}\otimes\omega_1)\oplus(1\otimes3\omega_1))$.\\
We have $G_v^\circ\cong\G_2$.
% 36
\item $(\GL_1^3\times\Spin_{10}\times\SL_2,(\mathrm{halfspinrep}\otimes\omega_1)\oplus(1\otimes\omega_1)\oplus(1\otimes\omega_1))$.\\
We have $G_v^\circ\cong\GL_1\times\G_2$.
The module $(\GL_1^2\times\Spin_{10}\times\SL_2,(\mathrm{halfspinrep}\otimes\omega_1)\oplus(1\otimes\omega_1)\oplus(1\otimes\omega_1))$
is prehomogeneous with $G_v^\circ\cong\G_2$.
% 37
\item $(\GL_1^3\times\Spin_{10}\times\SL_2,(\mathrm{halfspinrep}\otimes\omega_1)\oplus(1\otimes2\omega_1)\oplus(1\otimes\omega_1))$.\\
We have $G_v^\circ\cong\G_2$.
% 38
\item $(\GL_1^4\times\Spin_{10}\times\SL_2,(\mathrm{halfspinrep}\otimes\omega_1)\oplus(1\otimes\omega_1)\oplus(1\otimes\omega_1)\oplus(1\otimes\omega_1))$.\\
We have $G_v^\circ\cong\G_2$.
% 39
\item $(\GL_1^2\times\Spin_{10}\times\SL_3,(\mathrm{halfspinrep}\otimes\omega_1)\oplus(1\otimes\omega_1^{(*)}))$.\\
We have $G_v^\circ\cong\SL_2\times\SO_2$.
% 40
\item $(\GL_1^2\times\Spin_{10}\times\SL_{14},(\mathrm{halfspinrep}\otimes\omega_1)\oplus(1\otimes\omega_1^*))$.\\
We have $G_v^\circ\cong\SL_2\times\SO_2$.
% 41
\item $(\GL_1^2\times\Spin_{10}\times\SL_{15},(\mathrm{halfspinrep}\otimes\omega_1)\oplus(1\otimes\omega_1^*))$.\\
We have $G_v^\circ\cong\GL_1\times\SL_4$.
The module $(\GL_1\times\Spin_{10}\times\SL_{15},(\mathrm{halfspinrep}\otimes\omega_1)\oplus(1\otimes\omega_1^{(*)}))$ is prehomogeneous
with $G_v^\circ\cong\SL_4$.
% 42
\item $(\GL_1^2\times\Spin_{10}\times\SL_2,(\mathrm{vecrep}\otimes\omega_1)\oplus(\mathrm{halfspinrep}\otimes 1))$.\\
We have $G_v^\circ\cong\G_2$.
% 43
\item $(\GL_1^2\times\Spin_{10}\times\SL_3,(\mathrm{halfspinrep}\otimes\omega_1)\oplus(\mathrm{vecrep}\otimes1))$.\\
We have $G_v^\circ\cong\SL_3\times\SO_2$.
% 44
\item $(\GL_1^2\times\Spin_{10}\times\SL_4,(\mathrm{halfspinrep}\otimes\omega_1)\oplus(\mathrm{vecrep}\otimes1))$.\\
We have $G_v^\circ\cong\SL_2\times\SL_2$.
% 45
\item $(\GL_1^2\times\G_2\times\SL_2,(\omega_2\otimes\omega_1)\oplus(1\otimes\omega_1))$.\\
We have $G_v^\circ\cong\SL_2$.
% 46
\item $(\GL_1^2\times\G_2\times\SL_6,(\omega_2\otimes\omega_1)\oplus(1\otimes\omega_1^*))$.\\
We have $G_v^\circ\cong\SL_2$.
\end{enumerate}
% KI II
\item[{\rm \textbf{KI II}}] Non-regular two-simple prehomogeneous modules of type I.
\begin{enumerate}
\item 
% 1
\begin{enumerate}
% 1a
\item $(\GL_1^2\times\SL_{2n+1}\times\SL_2,
(\omega_2\otimes\omega_1)\oplus(1\otimes \omega_1))$ for $n\geq 2$.\\
We have $G_v^\circ\cong \GL_1^2\cdot \G^+_1$.
% 1b
\item $(\GL_1^2\times\SL_{2n+1}\times\SL_2,
(\omega_2\otimes\omega_1)\oplus(1\otimes 2\omega_1))$ for $n\geq 2$.\\
We have $G_v^\circ\cong (\GL_1\times\SO_2)\cdot\G^+_{2n}$.
% 1c
\item $(\GL_1^2\times\SL_{2n+1}\times\SL_2,
(\omega_2\otimes\omega_1)\oplus(1\otimes 3\omega_1))$ for $n\geq 2$.\\
We have $G_v^\circ\cong \GL_1\cdot\G^+_n$.
\end{enumerate}
% 2
\item 
\begin{enumerate}
% 2a
\item $(\GL_1^3\times\SL_{2n+1}\times\SL_2,
(\omega_2\otimes\omega_1)\oplus(1\otimes \omega_1)\oplus(1\otimes \omega_1))$ for $n\geq 2$.\\
We have $G_v^\circ\cong \GL_1^2\cdot\G^+_{2n}$.
% 2b
\item $(\GL_1^3\times\SL_{2n+1}\times\SL_2,
(\omega_2\otimes\omega_1)\oplus(1\otimes \omega_1)\oplus(1\otimes2\omega_1))$ for $n\geq 2$.\\
We have $G_v^\circ\cong \GL_1\cdot \G^+_{2n}$.
\end{enumerate}
% 3
\item $(\GL_1^4\times\SL_{2n+1}\times\SL_2,(\omega_2\otimes\omega_1)\oplus(1\otimes\omega_1)\oplus(1\otimes\omega_1)\oplus(1\otimes\omega_1))$ for $n\geq 2$.\\
We have $G_v^\circ\cong\GL_1\cdot\G^+_{2n}$.
% 4
\item $(\GL_1^2\times\SL_4\times\SL_2,(\omega_2\otimes\omega_1)\oplus(\omega_1\otimes 1))$.\\
We have $G_v^\circ\cong(\GL_1\times\SO_2)\cdot\Un_2$.
% 5
\item $(\GL_1^3\times\SL_4\times\SL_2,(\omega_2\otimes\omega_1)\oplus(\omega_1\otimes 1)\oplus(1\otimes\omega_1))$.\\
We have $G_v^\circ\cong\GL_1\cdot\Un_2$.
% 6
\item $(\GL_1^3\times\SL_4\times\SL_2,(\omega_2\otimes\omega_1)\oplus(\omega_1\otimes 1)\oplus(1\otimes\omega_1^*))$.\\
We have $G_v^\circ\cong\GL_1\cdot\Un_2$.
% 7
\item
\begin{enumerate}
% 7a
\item $(\GL_1^2\times\SL_5\times\SL_2,(\omega_2\otimes\omega_1)\oplus(\omega_1\otimes 1))$.\\
We have $G_v^\circ\cong(\GL_1\times\SO_2)\cdot\Un_2$.
% 7b
\item $(\GL_1^2\times\SL_5\times\SL_2,(\omega_2\otimes\omega_1)\oplus(\omega_1^*\otimes 1))$.\\
We have $G_v^\circ\cong\GL_1^2\cdot\Un_2$.
\end{enumerate}
% 8
\item
\begin{enumerate}
% 8a
\item $(\GL_1^3\times\SL_5\times\SL_2,(\omega_2\otimes\omega_1)\oplus(\omega_1\otimes 1)\oplus(1\otimes\omega_1))$.\\
We have $G_v^\circ\cong\GL_1\cdot\Un_2$.
% 8b
\item $(\GL_1^3\times\SL_5\times\SL_2,(\omega_2\otimes\omega_1)\oplus(\omega_1^*\otimes 1)\oplus(1\otimes\omega_1))$.\\
We have $G_v^\circ\cong\GL_1^2\cdot\Un_2$.
\end{enumerate}
% 9
\item $(\GL_1^3\times\SL_5\times\SL_9,(\omega_2\otimes\omega_1)\oplus(\omega_1^{(*)}\otimes 1)\oplus(1\otimes\omega_1^*))$.\\
We have $G_v^\circ\cong\GL_1\cdot\Un_2$.
% 10
\item $(\GL_1^3\times\SL_5\times\SL_2,(\omega_2\otimes\omega_1)\oplus(\omega_1^*\otimes 1)\oplus(1\otimes 2\omega_1))$.\\
We have $G_v^\circ\cong\GL_1\cdot\Un_2$.
% 11
\item $(\GL_1^4\times\SL_5\times\SL_2,(\omega_2\otimes\omega_1)\oplus(\omega_1^*\otimes 1)\oplus(1\otimes \omega_1)\oplus(1\otimes\omega_1))$.\\
We have $G_v^\circ\cong\GL_1\cdot\Un_2$.
% 12
\item $(\GL_1^2\times\SL_6\times\SL_2,(\omega_2\otimes\omega_1)\oplus(\omega_1^{(*)}\otimes 1))$.\\
We have $G_v^\circ\cong\GL_1\cdot\Un_3$.
% 13
\item
\begin{enumerate}
% 13a
\item $(\GL_1^2\times\SL_7\times\SL_2,(\omega_2\otimes\omega_1)\oplus(\omega_1\otimes1))$.\\
We have $G_v^\circ\cong\GL_1\cdot\Un_3$.
% 13b
\item $(\GL_1^2\times\SL_7\times\SL_2,(\omega_2\otimes\omega_1)\oplus(\omega_1^*\otimes1))$.\\
We have $G_v^\circ\cong(\GL_1\times\SO_2)\cdot\Un_2$.
\end{enumerate}
% 14
\item $(\GL_1^3\times\SL_7\times\SL_2,(\omega_2\otimes\omega_1)\oplus(\omega_1^*\otimes1)\oplus(1\otimes\omega_1)$.\\
We have $G_v^\circ\cong\GL_1\cdot\Un_2$.
% 15
\item $(\GL_1^2\times\SL_9\times\SL_2,(\omega_2\otimes\omega_1)\oplus(\omega_1^*\otimes1))$.\\
We have $G_v^\circ\cong\GL_1\cdot\Un_3$.
% 16
\item
\begin{enumerate}
% 16a
\item $(\GL_1^2\times\Sp_n\times\SL_{2m},(\omega_1\otimes\omega_1)\oplus(\omega_1\otimes 1))$ for $n>m\geq 1$.\\
We have $G_v^\circ\cong(\GL_1\times\Sp_{m-1}\times\Sp_{n-m-1})\cdot\Un_{2n-2}$.
% 16b
\item $(\GL_1^2\times\Sp_n\times\SL_{2m},(\omega_1\otimes\omega_1)\oplus(1\otimes\omega_1^{(*)}))$ for $n>m\geq 1$.\\
We have $G_v^\circ\cong(\GL_1\times\Sp_{m-1}\times\Sp_{n-m-1})\cdot\Un_{2n-2}$.

\end{enumerate}
% 17
\item $(\GL_1^3\times\Sp_n\times\SL_{2m},(\omega_1\otimes\omega_1)\oplus(\omega_1\otimes1)\oplus(1\otimes\omega_1^{(*)}))$ for $n>m\geq 1$.\\
We have $G_v^\circ\cong(\GL_1\times\Sp_{m-1}\times\Sp_{n-m-1})\cdot\Un_{2n-2m-2}$.
% 18
\item
\begin{enumerate}
% 18ax ay az
\item $(\GL_1^4\times\Sp_n\times\SL_{2m},(\omega_1\otimes\omega_1)\oplus(\omega_1\otimes 1)\oplus(1\otimes\omega_1^{(*)})\oplus(1\otimes\omega_1^{(*)}))$ for $n>m\geq 2$.\\
We have $G_v^\circ\cong(\Sp_{m-2}\times\Sp_{n-m-1})\cdot\Un_{2n-4}$.
% 18aw
\item $(\GL_1^4\times\Sp_n\times\SL_2,(\omega_1\otimes\omega_1)\oplus(\omega_1\otimes 1)\oplus(1\otimes\omega_1)\oplus(1\otimes\omega_1))$ for $n\geq 2$.\\
We have $G_v^\circ\cong\Sp_{n-2}\cdot\Un_{2n-3}$.
% 18bx by bz bw
\item $(\GL_1^4\times\Sp_n\times\SL_{2m},(\omega_1\otimes\omega_1)\oplus(1\otimes\omega_1^{(*)})\oplus(1\otimes\omega_1^{(*)})\oplus(1\otimes\omega_1^{(*)}))$ for $n>m\geq 2$.\\
We have $G_v^\circ\cong(\Sp_{m-2}\times\Sp_{n-m-1})\cdot\Un_{2n-4}$.
\end{enumerate}
% 19
\item $(\GL_1^3\times\Sp_n\times\SL_2,
(\omega_1\otimes\omega_1)\oplus(\omega_1\otimes 1)\oplus(1\otimes2\omega_1))$ for $n\geq 2$.\\
We have $G_v^\circ\cong\Sp_{n-2}\cdot\Un_{2n-3}$.
% 20
\item
\begin{enumerate}
% 20a
\item $(\GL_1^2\times\Sp_n\times\SL_{2m+1},
(\omega_1\otimes\omega_1)\oplus(1\otimes\omega_1))$ for $n>m\geq 1$.\\
We have $G_v^\circ\cong(\GL_1\times\Sp_m\times\Sp_{n-m-1})\cdot\Un_{2n-2m-1}$.
% 20b
\item $(\GL_1^2\times\Sp_n\times\SL_{2m+1},
(\omega_1\otimes\omega_1)\oplus(1\otimes\omega_1^*))$ for $n>m\geq 1$.\\
We have $G_v^\circ\cong(\GL_1^2\times\Sp_m\times\Sp_{n-m-1})\cdot\Un_{2n-2m-3}$.
% 20b'
\item $(\GL_1^2\times\Sp_{m+1}\times\SL_{2m+1},
(\omega_1\otimes\omega_1)\oplus(1\otimes\omega_1))$ for $m\geq 1$.\\
We have $G_v^\circ\cong(\GL_1^2\times\Sp_{m-1})\cdot\Un_{4m-1}$.
% 20c
\item $(\GL_1^2\times\Sp_n\times\SL_{2m+1},
(\omega_1\otimes\omega_1)\oplus(1\otimes\omega_2))$ for $n>m+1\geq 2$.\\
We have $G_v^\circ\cong(\GL_1\times\SO_2^m\times\Sp_{n-m-1})\cdot\Un_{2n-2m-1}$.
\end{enumerate}
% 21
\item
\begin{enumerate}
% 21a
\item $(\GL_1^3\times\Sp_n\times\SL_{2m+1},
(\omega_1\otimes\omega_1)\oplus(\omega_1\otimes1)\oplus(\omega_1\otimes1))$ for $n>m+1$.\\
We have $G_v^\circ\cong(\Sp_{m-1}\times\Sp_{n-m-2})\cdot\Un_{2n-4}$.
% 21b
\item $(\GL_1^3\times\Sp_n\times\SL_{2m+1},
(\omega_1\otimes\omega_1)\oplus(\omega_1\otimes1)\oplus(1\otimes\omega_1^{(*)}))$ for $n>m+1$.\\
We have $G_v^\circ\cong(\Sp_{m-1}\times\Sp_{n-m-2})\cdot\Un_{2n-4}$.
% 21c d e
\item $(\GL_1^3\times\Sp_n\times\SL_{2m+1},
(\omega_1\otimes\omega_1)\oplus(1\otimes\omega_1^{(*)})\oplus(1\otimes\omega_1^{(*)}))$ for $n>m+1$.\\
We have $G_v^\circ\cong(\Sp_{m-1}\times\Sp_{n-m-2})\cdot\Un_{2n-4}$.
\end{enumerate}
% 22
\item
\begin{enumerate}
% 22a b
\item $(\GL_1^4\times\Sp_n\times\SL_{2m+1},
(\omega_1\otimes\omega_1)\oplus(1\otimes\omega_1)\oplus
(1\otimes(\omega_1\oplus\omega_1)^{(*)}))$ for $n>m\geq 1$.\\
We have $G_v^\circ\cong(\Sp_{m-1}\times\Sp_{n-m-1})\cdot\Un_{2n-2m-1}$.
% 22c
\item $(\GL_1^4\times\Sp_n\times\SL_{2m+1},
(\omega_1\otimes\omega_1)\oplus
(1\otimes\omega_1^*)\oplus(1\otimes\omega_1^*)\oplus(1\otimes\omega_1^*))$ for $n>m\geq 1$.\\
We have $G_v^\circ\cong(\Sp_{m-2}\times\Sp_{n-m-1})\cdot\Un_{2n+2m-7}$.
% 22d
\item $(\GL_1^4\times\Sp_n\times\SL_3,
(\omega_1\otimes\omega_1)\oplus
(1\otimes\omega_1^*)\oplus(1\otimes\omega_1^*)\oplus(1\otimes\omega_1^*))$ for $n\geq 2$.\\
We have $G_v^\circ\cong\Sp_{n-2}\cdot\Un_{2n-3}$.
\end{enumerate}
% 23
\item $(\GL_1^2\times\Sp_n\times \SL_3,(\omega_1\otimes\omega_1)\oplus(1\otimes 2\omega_1))$ for $n\geq 2$.\\
We have $G_v^\circ\cong(\SO_2\times\Sp_{n-2})\cdot\Un_{2n-3}$.
% 24
\item $(\GL_1^3\times\Sp_n\times \SL_5,(\omega_1\otimes\omega_1)\oplus(1\otimes \omega_2)\oplus(1\otimes\omega_1^*))$ for $n\geq 3$.\\
We have $G_v^\circ\cong(\GL_1\times\Sp_{n-3})\cdot\Un_{2n-4}$.
% 25
\item $(\GL_1^2\times\Sp_n\times \SL_2,(\omega_1\otimes2\omega_1)\oplus(1\otimes \omega_1))$ for $n\geq 2$.\\
We have $G_v^\circ\cong\Sp_{n-2}\cdot\Un_{2n-3}$.
% 26
\item $(\GL_1^2\times\Spin_{10}\times \SL_2,(\mathrm{halfspinrep}\otimes\omega_1)\oplus(1\otimes \omega_1))$.\\
We have $G_v^\circ\cong(\GL_1\times\G_2)\cdot\G^+_1$.
\end{enumerate}
\end{itemize}

\subsection{Two-simple prehomogeneous modules of type II}

In this chapter we give a classification of the two-simple prehomogeneous modules of type II,
i.e. modules of the form
\begin{align*}
\bigl(&\GL_1^l\times G_1 \times G_2,\\
&(\rho_1 \otimes \tilde{\rho}_1)\oplus\ldots\oplus(\rho_k\otimes\tilde{\rho}_k)\
\oplus\
(\sigma_1\otimes 1)\oplus\ldots\oplus(\sigma_s\otimes 1)\
\oplus\
(1\otimes\tau_1)\oplus\ldots\oplus(1\otimes\tau_t),\\
&V_1\oplus\ldots \oplus V_l \bigr),
\end{align*}
where all of the modules
$(\GL_1\times G_1\times G_2, \rho_i\otimes\tilde{\rho}_i)$
are trivial prehomogeneous modules (see Kimura \cite{kimura}).
Note that we consider non-irreducible modules only.
These were classified by Kimura et al. \cite{kimuraII}, thus we
shall refer to them as KII $n$, where $n$ is the number of the module
in \S 5 of \cite{kimuraII}. 
Unfortunately, it is not always obvious from the classification
in which cases a module would be prehomogeneous even with fewer
than $l$ scalar multiplications.
\index{type II}
\index{KII}
\index{prehomogeneous module!2-simple}
\index{prehomogeneous module!classification}

Any indecomposable two-simple prehomogeneous module of type II
is  equivalent to one of the following:
\begin{itemize}
% KII I
\item[{\rm \textbf{KII I}}] Two-simple prehomogeneous modules of type II obtained
directly from any given simple module
$(\GL_1^l\times G,\rho_1\oplus\ldots\oplus\rho_l)$.
\begin{enumerate}
% 1
\item For any representation $\sigma_1\oplus\ldots\oplus\sigma_s$
of $G$ and $n\geq \sum_{i=1}^s\dim(\sigma_i)$:
\begin{align*}
\bigl(&\GL_1^{l+s}\times G\times \SL_n,\\
&(\sigma_1\otimes\omega_1)\oplus\ldots\oplus(\sigma_s\otimes\omega_1)\
\oplus\ (\rho_1\otimes 1)\oplus\ldots\oplus(\rho_l\otimes 1)\bigr).
\end{align*}
% 2
\item For $t\geq 0$, $1\leq k\leq l$ and
$n=t-1+\sum_{i=1}^k \dim(\rho_i)$:
\begin{align*}
\bigl(&\GL_1^{l+t}\times G\times\SL_n,\\
&(\rho_1\otimes\omega_1)\oplus\ldots\oplus(\rho_k\otimes\omega_1)\
\oplus\ (\rho_{k+1}^*\otimes 1)\oplus\ldots\oplus(\rho_l^*\otimes 1)\
\oplus\ (1\otimes\omega_1^{\oplus t}) \bigr).
\end{align*}
% 3
\item For $t\geq 1$, $1\leq k\leq l$ and 
$n\geq t-1+\sum_{i=1}^k \dim(\rho_i)$:
\begin{align*}
\bigl( &\GL_1^{l+t}\times G\times \SL_n,\\
&(\rho_1\otimes\omega_1)\oplus\ldots\oplus(\rho_k\otimes\omega_1)\
\oplus\ (\rho_{k+1}\otimes 1)\oplus\ldots\oplus(\rho_l\otimes 1)\
\oplus\ (1\otimes\omega_1^{\oplus t-1})\oplus(1\otimes\omega_1^*) 
\bigr).
\end{align*}
\end{enumerate}
% KII II
\item[{\rm \textbf{KII II}}] Two-simple prehomogeneous modules of type II of
the form
\begin{align*}
\bigl(&\GL_1^{k+s+t}\times G \times \SL_n,\\
&(\rho_1 \otimes \omega_1)\oplus\ldots\oplus(\rho_k\otimes\omega_1)\
\oplus\
(\sigma_1\otimes 1)\oplus\ldots\oplus(\sigma_s\otimes 1)\
\oplus\
(1\otimes\tau_1)\oplus\ldots\oplus(1\otimes\tau_t)
\bigr),
\end{align*}
with $2\leq\dim(\rho_i)\leq n$ for all $i$ and at least one
$\tau_j\neq\omega_1^{(*)}$.
\begin{enumerate}\setcounter{enumi}{3}
% 4
\item $G=\SL_m$ with $2\leq m < n$.
\begin{itemize}
\item[4-i] 
\begin{enumerate}
\item $(\omega_1\otimes\omega_1)\oplus(1\otimes2\omega_1^{(*)})$.
\item $(\omega_1\otimes\omega_1)\oplus(1\otimes2\omega_1^{(*)})\oplus(\omega_1^{(*)}\otimes 1)$.
\item $(\omega_1\otimes\omega_1)\oplus(1\otimes\omega_2^{(*)})$.
\item $(\omega_1\otimes\omega_1)\oplus(1\otimes\omega_2^{(*)})\oplus(\omega_1^{(*)}\otimes 1)$.
\item $(\omega_1\otimes\omega_1)\oplus(1\otimes\omega_2^{(*)})\oplus(\omega_1^{(*)}\otimes 1)\oplus(\omega_1^{(*)}\otimes 1)$.
\item $(\omega_1\otimes\omega_1)\oplus(1\otimes\omega_2^{(*)})\oplus(1\otimes\omega_1^{(*)})$.
\item $(\omega_1\otimes\omega_1)\oplus(1\otimes\omega_2^{(*)})\oplus(\omega_1^{(*)}\otimes 1)\oplus(1\otimes\omega_1^{(*)})$.
\end{enumerate}
\item[4-ii] $n$ even.
\begin{enumerate}
\item $(\omega_1\otimes\omega_1)\oplus(1\otimes\omega_2^{(*)})\oplus((\omega_1\oplus\omega_1)^{(*)}\otimes 1)\oplus(1\otimes\omega_1^{(*)})$.
\item $(\omega_1\otimes\omega_1)\oplus(1\otimes\omega_2^{(*)})\oplus(\omega_1\otimes1)\oplus(\omega_1^*\otimes 1)\oplus(\omega_1^*\otimes 1)$.
\item $(\omega_1\otimes\omega_1)\oplus(1\otimes\omega_2^{(*)})\oplus((\omega_1\oplus\omega_1\oplus\omega_1)^{(*)}\otimes 1)$.
\item $(\omega_1\otimes\omega_1)\oplus(1\otimes\omega_2^{(*)})\oplus(\omega_1\otimes 1)\oplus(\omega_1\otimes1)\oplus(\omega_1^*\otimes1)$, $m$ even.
\item $(\omega_1\otimes\omega_1)\oplus(1\otimes\omega_2^{(*)})\oplus(\omega_1\otimes1)\oplus(\omega_1^*\otimes1)\oplus(1\otimes\omega_1^{(*)})$, $m$ even.
\item $(\omega_1\otimes\omega_1)\oplus(1\otimes\omega_2^{(*)})\oplus(\omega_2\otimes1)$, $m$ odd.
\item $(\omega_1\otimes\omega_1)\oplus(1\otimes\omega_2^{(*)})\oplus(1\otimes\omega_1^{(*)})\oplus(1\otimes\omega_1^{(*)})$, $m$ odd.
\end{enumerate}
\item[4-iii] $n$ odd.
\begin{enumerate}
\item $(\omega_1\otimes\omega_1)\oplus(1\otimes\omega_2)\oplus(\omega_1\otimes1)\oplus(\omega_1\otimes1)\oplus(1\otimes\omega_1^*)$, $m\geq 3$.
\item $(\omega_1\otimes\omega_1)\oplus(1\otimes\omega_2)\oplus(\omega_2^*\otimes 1)$.
\item $(\omega_1\otimes\omega_1)\oplus(1\otimes\omega_2)\oplus((\omega_1\oplus\omega_1\oplus\omega_1)^{(*)}\otimes 1)$, $m\geq 3$.
\item $(\omega_1\otimes\omega_1)\oplus(1\otimes\omega_2)\oplus((\omega_1\oplus\omega_1)^{(*)}\otimes1)\oplus(1\otimes\omega_1)$.
\item $(\omega_1\otimes\omega_1)\oplus(1\otimes\omega_2)\oplus(\omega_2\otimes1)$, $m$ even.
\item $(\omega_1\otimes\omega_1)\oplus(1\otimes\omega_2)\oplus(\omega_1\otimes1)\oplus(\omega_1^*\otimes1)\oplus(1\otimes\omega_1)$, $m$ even.
\item $(\omega_1\otimes\omega_1)\oplus(1\otimes\omega_2)\oplus(1\otimes\omega_1)\oplus(1\otimes\omega_1^*)$, $m$ even.
\item $(\omega_1\otimes\omega_1)\oplus(1\otimes\omega_2)\oplus(1\otimes\omega_1^*)\oplus(1\otimes\omega_1^*)$, $m$ even.
\item $(\omega_1\otimes\omega_1)\oplus(1\otimes\omega_2)\oplus(\omega_1\otimes1)\oplus(\omega_1\otimes1)\oplus(\omega_1^*\otimes1)$, $m$ odd.
\item $(\omega_1\otimes\omega_1)\oplus(1\otimes\omega_2)\oplus(1\otimes\omega_1)\oplus(1\otimes\omega_1)$, $m$ odd.
\item $(\omega_1\otimes\omega_1)\oplus(1\otimes\omega_2)\oplus(\omega_1\otimes1)\oplus(\omega_1^*\otimes1)\oplus(1\otimes\omega_1^*)$, $m$ odd.
\item $(\omega_1\otimes\omega_1)\oplus(1\otimes\omega_2^*)\oplus(\omega_1\otimes1)\oplus(\omega_1^*\otimes1)\oplus(\omega_1^*\otimes1)$.
\item $(\omega_1\otimes\omega_1)\oplus(1\otimes\omega_2^*)\oplus((\omega_1\oplus\omega_1\oplus\omega_1)^{(*)}\otimes1)$.
\item $(\omega_1\otimes\omega_1)\oplus(1\otimes\omega_2^*)\oplus((\omega_1\oplus\omega_1)^{(*)}\otimes1)\oplus(1\otimes\omega_1^{(*)})$.
\item $(\omega_1\otimes\omega_1)\oplus(1\otimes\omega_2^*)\oplus(\omega_1\otimes1)\oplus(\omega_1\otimes1)\oplus(\omega_1^*\otimes1)$, $m$ even.
\item $(\omega_1\otimes\omega_1)\oplus(1\otimes\omega_2^*)\oplus(1\otimes\omega_1^*)\oplus(1\otimes\omega_1^*)$, $m$ even.
\item $(\omega_1\otimes\omega_1)\oplus(1\otimes\omega_2^*)\oplus(\omega_1\otimes1)\oplus(\omega_1^*\otimes1)\oplus(1\otimes\omega_1^{(*)})$, $m$ even.
\item $(\omega_1\otimes\omega_1)\oplus(1\otimes\omega_2^*)\oplus(\omega_2\otimes1)$, $m$ odd.
\item $(\omega_1\otimes\omega_1)\oplus(1\otimes\omega_2^*)\oplus(1\otimes\omega_1)\oplus(1\otimes\omega_1)$, $m$ odd.
\item $(\omega_1\otimes\omega_1)\oplus(1\otimes\omega_2^*)\oplus(1\otimes\omega_1)\oplus(1\otimes\omega_1^*)$, $m$ odd.
\end{enumerate}
\end{itemize}
% 5
\item $G=\SL_2$, $n > 2$.
\begin{itemize}
\item[5-i]
\begin{enumerate}
\item $(2\omega_1\otimes\omega_1)\oplus(1\otimes\omega_2^{(*)})$.
\item $(2\omega_1\otimes\omega_1)\oplus(1\otimes\omega_2^{(*)})\oplus(\omega_1\otimes1)$.
\end{enumerate}
\item[5-ii]
\begin{enumerate}
\item $(\omega_1\otimes\omega_1)\oplus(1\otimes\omega_2)\oplus(2\omega_1\otimes1)$.
\item $(\omega_1\otimes\omega_1)\oplus(1\otimes\omega_2)\oplus(2\omega_1\otimes1)\oplus(1\otimes\omega_1)$.
\item $(\omega_1\otimes\omega_1)\oplus(1\otimes\omega_2)\oplus(3\omega_1\otimes1)$, $n$ even.
\item $(\omega_1\otimes\omega_1)\oplus(1\otimes\omega_2)\oplus(2\omega_1\otimes1)\oplus(\omega_1\otimes1)$, $n$ even.
\item $(\omega_1\otimes\omega_1)\oplus(1\otimes\omega_2)\oplus(2\omega_1\otimes1)\oplus(1\otimes\omega_1^*)$, $n$ even.
\end{enumerate}
\item[5-iii]
\begin{enumerate}
\item $(\omega_1\otimes\omega_1)\oplus(1\otimes\omega_2^*)$.
\item $(\omega_1\otimes\omega_1)\oplus(1\otimes\omega_2^*)\oplus(2\omega_1\otimes1)$.
\item $(\omega_1\otimes\omega_1)\oplus(1\otimes\omega_2^*)\oplus(3\omega_1\otimes1)$.
\item $(\omega_1\otimes\omega_1)\oplus(1\otimes\omega_2^*)\oplus(2\omega_1\otimes1)\oplus(\omega_1\otimes1)$.
\item $(\omega_1\otimes\omega_1)\oplus(1\otimes\omega_2^*)\oplus(2\omega_1\otimes1)\oplus(1\otimes\omega_1^{(*)})$.
\end{enumerate}
\item[5-iv] $n=5$.
\begin{enumerate}
\item $(\omega_1\otimes\omega_1)\oplus(1\otimes\omega_2^*)\oplus(1\otimes\omega_2^*)$.
\end{enumerate}
\item[5-v] $n=6$.
\begin{enumerate}
\item $(\omega_1\otimes\omega_1)\oplus(1\otimes\omega_3)$.
\item $(\omega_1\otimes\omega_1)\oplus(1\otimes\omega_3)\oplus(1\otimes\omega_1^*)$.
\item $(\omega_1\otimes\omega_1)\oplus(1\otimes\omega_3)\oplus(\omega_1\otimes1)$.
\item $(\omega_1\otimes\omega_1)\oplus(1\otimes\omega_3)\oplus(2\omega_1\otimes1)$.
\item $(\omega_1\otimes\omega_1)\oplus(1\otimes\omega_3)\oplus(3\omega_1\otimes1)$.
\item $(\omega_1\otimes\omega_1)\oplus(1\otimes\omega_3)\oplus(\omega_1\otimes1)\oplus(\omega_1\otimes1)$.
\item $(\omega_1\otimes\omega_1)\oplus(1\otimes\omega_3)\oplus(2\omega_1\otimes1)\oplus(\omega_1\otimes1)$.
\item $(\omega_1\otimes\omega_1)\oplus(1\otimes\omega_3)\oplus(\omega_1\otimes1)\oplus(\omega_1\otimes1)\oplus(\omega_1\otimes1)$.
\end{enumerate}
\item[5-vi] $n=7$.
\begin{enumerate}
\item $(\omega_1\otimes\omega_1)\oplus(1\otimes\omega_3^{(*)})$.
\item $(\omega_1\otimes\omega_1)\oplus(1\otimes\omega_3^{(*)})\oplus(\omega_1\otimes1)$.
\end{enumerate}
\end{itemize}
% 6
\item $G=\SL_3$, $n>3$.
\begin{enumerate}
\item $(\omega_1\otimes\omega_1)\oplus(1\otimes\omega_2^{(*)})\oplus(2\omega_1^{(*)}\otimes1)$.
\item $(\omega_1\otimes\omega_1)\oplus(1\otimes\omega_2)\oplus(1\otimes\omega_2)$, $n=5$.
\end{enumerate}
% 7
\item $G=\SL_4$, $n>4$.
\begin{itemize}
\item[7-i] $n$ odd.
\begin{enumerate}
\item $(\omega_1\otimes\omega_1)\oplus(1\otimes\omega_2)\oplus(2\omega_1\otimes1)$.
\item $(\omega_1\otimes\omega_1)\oplus(1\otimes\omega_2)\oplus(\omega_2\otimes1)\oplus(\omega_1^{(*)}\otimes1)$.
\item $(\omega_1\otimes\omega_1)\oplus(1\otimes\omega_2)\oplus(\omega_2\otimes1)\oplus(1\otimes\omega_1^*)$.
\end{enumerate}
\item[7-ii] $n=5$.
\begin{enumerate}
\item $(\omega_1\otimes\omega_1)\oplus(1\otimes\omega_2)\oplus(1\otimes\omega_2)$.
\end{enumerate}
\item[7-iii] $n=6$.
\begin{enumerate}
\item $(\omega_1\otimes\omega_1)\oplus(1\otimes\omega_3)$.
\item $(\omega_1\otimes\omega_1)\oplus(1\otimes\omega_3)\oplus(\omega_1^*\otimes1)$.
\item $(\omega_1\otimes\omega_1)\oplus(1\otimes\omega_3)\oplus(\omega_2^{(*)}\otimes1)$.
\item $(\omega_1\otimes\omega_1)\oplus(1\otimes\omega_3)\oplus(\omega_1^*\otimes1)\oplus(\omega_1^*\otimes1)$.
\item $(\omega_1\otimes\omega_1)\oplus(1\otimes\omega_3)\oplus(1\otimes\omega_1)$.
\end{enumerate}
\end{itemize}
% 8
\item $G=\SL_5$, $n>5$.
\begin{itemize}
\item[8-i] $n$ even.
\begin{enumerate}
\item $(\omega_1\otimes\omega_1)\oplus(1\otimes\omega_2^{(*)})\oplus(\omega_2\otimes1)\oplus(\omega_1^*\otimes1)$.
\end{enumerate}
\item[8-ii] $n$ odd.
\begin{enumerate}
\item $(\omega_1\otimes\omega_1)\oplus(1\otimes\omega_2^*)\oplus(\omega_2\otimes1)\oplus(\omega_1^*\otimes1)$.
\item $(\omega_1\otimes\omega_1)\oplus(1\otimes\omega_2^*)\oplus(1\otimes\omega_2)\oplus(\omega_2^*\otimes1)\oplus(\omega_1\otimes1)$.
\item $(\omega_1\otimes\omega_1)\oplus(1\otimes\omega_2^*)\oplus(1\otimes\omega_2)\oplus(\omega_2^*\otimes1)\oplus(1\otimes\omega_1^*)$.
\end{enumerate}
\item[8-iii] $n=6$.
\begin{enumerate}
\item $(\omega_1\otimes\omega_1)\oplus(1\otimes\omega_3)$.
\item $(\omega_1\otimes\omega_1)\oplus(1\otimes\omega_3)\oplus(\omega_1^{(*)}\otimes1)$.
\item $(\omega_1\otimes\omega_1)\oplus(1\otimes\omega_3)\oplus(\omega_2^*\otimes1)$.
\item $(\omega_1\otimes\omega_1)\oplus(1\otimes\omega_3)\oplus(\omega_1^*\otimes1)\oplus(\omega_1^{(*)}\otimes1)$.
\item $(\omega_1\otimes\omega_1)\oplus(1\otimes\omega_3)\oplus(1\otimes\omega_1^*)$.
\end{enumerate}
\item[8-iv] $n=7$.
\begin{enumerate}
\item $(\omega_1\otimes\omega_1)\oplus(1\otimes\omega_3^{(*)})$.
\end{enumerate}
\end{itemize}
% 9
\item $G=\SL_{2j}$, $n=2j+1$.
\begin{enumerate}
\item $(\omega_1\otimes\omega_1)\oplus(1\otimes\omega_2^*)\oplus(1\otimes\omega_1)\oplus(1\otimes\omega_1)$.
\item $(\omega_1\otimes\omega_1)\oplus(1\otimes\omega_3^{(*)})\oplus(\omega_1^{(*)}\otimes 1)$, $j=3$ (i.e. $n=7$).
\end{enumerate}
% 10
\item $G=\SL_n$.
\[
(\omega_1\otimes\omega_1)\
\oplus\ (\rho_1\otimes1)\oplus\ldots\oplus(\rho_k\otimes1)\
\oplus\ (1\otimes\rho_{k+1}^*)\oplus\ldots\oplus(1\otimes\rho_r^*),
\]
where $(\GL_1^r\times\SL_n,\rho_1\oplus\ldots\oplus\rho_r)$ is a
simple prehomogeneous module.
% 11
\item $G=\Sp_m$, $2m<n$.
\begin{itemize}
\item[11-i] $n$ odd.
\begin{enumerate}
\item $(\omega_1\otimes\omega_1)\oplus(1\otimes\omega_2)$.
\item $(\omega_1\otimes\omega_1)\oplus(1\otimes\omega_2)\oplus(\omega_1\otimes1)$, $m=2$.
\item $(\omega_1\otimes\omega_1)\oplus(1\otimes\omega_2)\oplus(\omega_2\otimes1)$, $m=2$.
\item $(\omega_1\otimes\omega_1)\oplus(1\otimes\omega_2)\oplus(1\otimes\omega_1^*)$, $m=2$.
\end{enumerate}
\item[11-ii] $n=6$, $m=2$.
\begin{enumerate}
\item $(\omega_1\otimes\omega_1)\oplus(1\otimes\omega_3)$.
\end{enumerate}
\end{itemize}
\end{enumerate}
% KII III
\item[{\rm \textbf{KII III}}] Two-simple prehomogeneous modules of type II of
the form
\begin{align*}
\bigl(&\GL_1^{k+s+t}\times G \times \SL_n,\\
&(\rho_1 \otimes \omega_1)\oplus\ldots\oplus(\rho_k\otimes\omega_1)\
\oplus\
(\sigma_1\otimes 1)\oplus\ldots\oplus(\sigma_s\otimes 1)\
\oplus\
(1\otimes\omega_1^{\oplus t})
\bigr),
\end{align*}
with $2\leq\dim(\rho_i)\leq n$ for all $i$ and 
\[(G,\ \rho_1,\ldots,\rho_k,\ \sigma_1,\ldots,\sigma_s)\neq
(\SL_m,\ \omega_1,\ldots,\omega_1,\ \omega_1^{(*)},\ldots,\omega_1^{(*)}).
\]
\begin{enumerate}\setcounter{enumi}{11}
% 12
\item $(\GL_1^2\times\SL_4\times\SL_8,(\omega_2\otimes\omega_1)\oplus(\omega_1\otimes\omega_1))$.
% 13
\item $G=\SL_m$.
\begin{itemize}
\item[13-i] $m<n$.
\begin{enumerate}
\item $(\omega_1\otimes\omega_1)\oplus(1\otimes\omega_1^*)\oplus(1\otimes\omega_1^*)\oplus(\omega_2^{(*)}\otimes1)$.
\item $(\omega_1\otimes\omega_1)\oplus(1\otimes\omega_1^*)\oplus(1\otimes\omega_1^*)\oplus(\omega_2^{(*)}\otimes1)\oplus(\omega_1\otimes1)$.
\item $(\omega_1\otimes\omega_1)\oplus(1\otimes\omega_1^*)\oplus(1\otimes\omega_1^*)\oplus(\omega_2^*\otimes1)\oplus(\omega_1^*\otimes1)$.
\item $(\omega_1\otimes\omega_1)\oplus(1\otimes\omega_1^*)\oplus(1\otimes\omega_1^*)\oplus(\omega_2^{(*)}\otimes1)\oplus(1\otimes\omega_1^*)$.
\item $(\omega_1\otimes\omega_1)\oplus(1\otimes\omega_1^*)\oplus(1\otimes\omega_1^*)\oplus(\omega_2^*\otimes1)\oplus(1\otimes\omega_1)$.
\item $(\omega_1\otimes\omega_1)\oplus(1\otimes\omega_1^*)\oplus(1\otimes\omega_1^*)\oplus(\omega_2\otimes1)\oplus(\omega_1^*\otimes1)$, $m$ even.
\item $(\omega_1\otimes\omega_1)\oplus(1\otimes\omega_1^*)\oplus(1\otimes\omega_1^*)\oplus(\omega_2\otimes1)\oplus(1\otimes\omega_1)$, $m$ even.
\item $(\omega_1\otimes\omega_1)\oplus(1\otimes\omega_1^*)\oplus(1\otimes\omega_1^*)\oplus(\omega_3\otimes1)$, $m=6$.
\item $(\omega_1\otimes\omega_1)\oplus(1\otimes\omega_1^*)\oplus(1\otimes\omega_1^*)\oplus(\omega_3\otimes1)\oplus(1\otimes\omega_1)$, $m=6$.
\end{enumerate}
\item[13-ii] $n=m+1$.
\begin{enumerate}
\item $(\omega_1\otimes\omega_1)\oplus(1\otimes\omega_1)\oplus(1\otimes\omega_1)\oplus(1\otimes\omega_1)\oplus(\omega_2\otimes1)$.
\item $(\omega_1\otimes\omega_1)\oplus(1\otimes\omega_1)\oplus(1\otimes\omega_1)\oplus(1\otimes\omega_1)\oplus(\omega_2^*\otimes1)$, $m$ even.
\item $(\omega_1\otimes\omega_1)\oplus(1\otimes\omega_1)\oplus(1\otimes\omega_1)\oplus(1\otimes\omega_1)\oplus(\omega_3\otimes1)$, $m=6$.
\end{enumerate}
\item[13-iii] $n\geq \frac{1}{2}m(m-1)$.
\begin{enumerate}
\item $(\omega_2\otimes\omega_1)\oplus(1\otimes\omega_1^*)\oplus(1\otimes\omega_1^*)$, $m$ odd.
\item $(\omega_2\otimes\omega_1)\oplus(1\otimes\omega_1)\oplus(1\otimes\omega_1^*)\oplus(1\otimes\omega_1^*)$, $m$ odd, $n>\frac{1}{2}m(m-1)$.
\item $(\omega_2\otimes\omega_1)\oplus(\omega_1^*\otimes1)\oplus(1\otimes\omega_1^*)\oplus(1\otimes\omega_1^*)$, $m=5$.
\item $(\omega_2\otimes\omega_1)\oplus(1\otimes\omega_1)\oplus(1\otimes\omega_1)\oplus(1\otimes\omega_1)$, $m=2j+1$, $n=2j^2+j+1$.
\item $(\omega_2\otimes\omega_1)\oplus(1\otimes\omega_1)\oplus(1\otimes\omega_1)$, $m=2j$, $n=2j^2+j$.
\item $(\omega_2\otimes\omega_1)\oplus(\omega_1\otimes1)\oplus(1\otimes\omega_1)\oplus(1\otimes\omega_1)$, $m=5$, $n=10$.
\end{enumerate}
\end{itemize}
% 14
\item $G=\Sp_m$, $n\geq 2m$.
\begin{enumerate}
\item $(\omega_1\otimes\omega_1)\oplus(1\otimes\omega_1^*)\oplus(1\otimes\omega_1^*)$.
\item $(\omega_1\otimes\omega_1)\oplus(1\otimes\omega_1^*)\oplus(1\otimes\omega_1^*)\oplus(\omega_1\otimes1)$.
\item $(\omega_1\otimes\omega_1)\oplus(1\otimes\omega_1^*)\oplus(1\otimes\omega_1^*)\oplus(1\otimes\omega_1^*)$.
\item $(\omega_1\otimes\omega_1)\oplus(1\otimes\omega_1^*)\oplus(1\otimes\omega_1^*)\oplus(1\otimes\omega_1)$, $n>2m$.
\item $(\omega_1\otimes\omega_1)\oplus(1\otimes\omega_1)\oplus(1\otimes\omega_1^*)\oplus(\omega_1\otimes1)$, $n=2m$.
\item $(\omega_1\otimes\omega_1)\oplus(1\otimes\omega_1)\oplus(1\otimes\omega_1^*)\oplus(1\otimes\omega_1^{(*)})$, $n=2m$.
\end{enumerate}
% 15
\item $G=\Spin_{10}$, $n\geq 16$.
\begin{enumerate}
\item $(\mathrm{halfspinrep}\otimes\omega_1)\oplus(1\otimes\omega_1^*)\oplus(1\otimes\omega_1^*)$.
\item $(\mathrm{halfspinrep}\otimes\omega_1)\oplus(1\otimes\omega_1)\oplus(1\otimes\omega_1^*)\oplus(1\otimes\omega_1^*)$, $n\geq 17$.
\item $(\mathrm{halfspinrep}\otimes\omega_1)\oplus(1\otimes\omega_1)\oplus(1\otimes\omega_1)\oplus(1\otimes\omega_1)$, $n=17$.
\item $(\mathrm{halfspinrep}\otimes\omega_1)\oplus(1\otimes\omega_1)\oplus(1\otimes\omega_1)$, $n=16$.
\end{enumerate}
\end{enumerate}
% KII IV
\item[{\rm \textbf{KII IV}}] Two-simple prehomogeneous modules of type II of
the form
\begin{align*}
\bigl(&\GL_1^{k+s_1+s_2+t_1+t_2}\times \SL_m \times \SL_n,\\
&(\omega_1\otimes 1)^{\oplus s_1}\oplus(\omega_1^*\otimes1)^{\oplus s_2}\
\oplus\ (\omega_1\otimes\omega_1)^{\oplus k}\
\oplus\ (1\otimes\omega_1)^{\oplus t_1}\oplus(1\otimes\omega_1^*)^{\oplus t_2}
\bigr),
\end{align*}
where $n\geq m\geq 2$ and $k\geq 1$.
\begin{enumerate}\setcounter{enumi}{15}
% 16
\item $n\geq km$.\label{thm_KII_IV}
\begin{itemize}
\item[16-i] $n=m$. Then $k=1$ and $1\leq(s_1+t_2)+(s_2+t_1)\leq n+1$,
where one of $s_1+t_2$ or $s_2+t_1$ is $0$ or $1$.
\item[16-ii] $n=km$, $k\geq2$.
\begin{enumerate}
\item $t_1=0$, $2\leq t_2\leq n$, $s_2=0$, $s_1+kt_2\leq m$.
\item $t_2=0$, $2\leq t_1\leq n$, $s_1=0$, $s_2+kt_1\leq m$.
\end{enumerate}
\item[16-iii] $n=km+1$. Then $t_1\geq 3$, $t_2=s_1=0$, $s_2+k(t_1-1)\leq m$.
\item[16-iv] $n\geq km+t_1$, $n>km$.
\begin{enumerate}
\item $k=1$, $t_1=0$, $2\leq t_2\leq n$ and $1\leq(s_1+t_2)+s_2\leq m+1$,
where $s_2$ is $0$ or $1$.
\item $k\geq 2$, $t_1=0$, $2\leq t_2\leq n$, $s_2=0$, $s_1+kt_2\leq m$.
\item $k\geq 1$, $t_1=1$, $2\leq t_2\leq n$, $s_2=0$, $s_1+kt_2\leq m$.
\end{enumerate}
\end{itemize}
% 17
\item $km > n$. 
These are the cases (17)-(25) in \S 5.4 of \cite{kimuraII}, but to
keep things simple we subsume them under the case KII IV-17 here.
See the following definition \ref{def_412} for the definition of
$T$, $\nu(k,m,n)$ and $(a_i)$.
Also, we write $b_i=\frac{a_i}{a_{i+1}}$.
\begin{itemize}
\item[17-i]
\begin{enumerate}
\item $t_2\geq 1$, $s_2=t_1=0$, $s_1+kt_2\leq m-b_j(n-t_2)$, where
$(k,m,n)\in T$ and $j=\nu(k,m,n)$.
\item $s_2=t_2=0$ and let $p=km+t_1-n(<m)$, $q=kp-m(<n)$,
$(k,p,m)\in T$ (resp.~$(k,q,p)\in T$) and
$j=\nu(k,p,m)$ (resp.~$j=(k,q,p)$).
\begin{enumerate}
\item $m\geq kp$, $s_1=0$ and $t_1\leq p+1$.
\item $m\geq kp$, $s_1=1$ and $k+t_1\leq p+1$.
\item $m\geq kp$, $2\leq s_1\leq m$ and $t_1+ks_1\leq p$.
\item $kp > m$, $s_1\geq 1$ and $t_1+ks_1\leq p-b_j(m-s_1)$.
\item $kp > m$, $p\geq kq$, $s_1=0$, $t_1=1$ and $k\leq q+1$.
\item $kp > m$, $p\geq kq$, $s_1=0$, $2\leq t_1\leq p$ and $kt_1\leq q$.
\item $kp > m$, $kq > p$, $s_1=0$, $t_1\geq 1$ and $kt_1\leq q-b_j(p-t_1)$.
\end{enumerate}
\item $t_2=0$, $s_2\geq 1$, $s_1=0$ and let $p=km+t_1-n(<m)$,
$q=kp+s_2-m(<p)$, $r=kq-p(<q)$, $(k,q,p)\in T$ (resp.~$(k,q,p)\in T$,
$(k,r,q)\in T$) and $j=\nu(k,q,p)$ (resp.~$j=\nu(k,r,q)$).
\begin{enumerate}
\item $m\geq kp+s_2$ and $t_1\leq p+1$.
\item $m=kp+s_2-1$, $t_1=0,1$ and $k+t_1\leq p+1$.
\item $m=kp+1$, $m\geq s_2\geq 3$, $t_1=0$ and $k(s_2-1)\leq p$.
\item $m=kp$, $m\geq s_2\geq 2$ and $ks_2\leq p$.
\item $kp>m$, $p\geq kq$, $t_1=0$ and $s_2\leq q+1$.
\item $kp>m$, $p\geq kq$, $t_1=1$ and $s_2+kt_1\leq q+1$.
\item $kp>m$, $p\geq kq$, $p\geq t_1\geq 2$ and $s_2+kt_1\leq q$.
\item $kp>m$, $kq>p$, $t_1\geq 1$ and $s_2+kt_1\leq q-b_j(p-t_1)$.
\item $kp>m$, $kq>p$, $q\geq kr$, $t_1=0$, $s_2=1$ and $k\leq r+1$.
\item $kp>m$, $kq>p$, $q\geq kr$, $t_1=0$, $q\geq s_2\geq 2$ and $ks_2\leq r$.
\item $kp>m$, $kq>p$, $kr>q$, $t_1=0$, $s_2\geq 1$ and $ks_2\leq r-b_j(q-s_2)$.
\end{enumerate}
\end{enumerate}
\item[17-ii]
\begin{enumerate}
\item $t_2=1$, $s_2=0$, $t_1\geq 1$ and let $p=km+t_1-n-1$,
$(k,p,m)\in T$ and $j=\nu(k,p,m)$.
\begin{enumerate}
\item $m\geq kp$ and $(t_1-1)+k(k+s_1)\leq p$.
\item $kp > m$ and $(t_1-1)+k(k+s_1)\leq p-b_j(m-k-s_1)$.
\end{enumerate}
\item $t_2=0$, $s_2\geq 1$, $s_1=1$ and let $p=km+t_1-n$,
$q=kp+s_2-m-1$, $(k,q,p)\in T$ and $j=\nu(k,q,p)$.
\begin{enumerate}
\item $kp> m$, $p\geq kq$, $(s_2-1)+k(k+t_1)\leq q$.
\item $kp> m$, $kq> n$, $(s_2-1)+k(k+t_1)\leq q-b_j(p-k-t_1)$.
\end{enumerate}
\end{enumerate}
\item[17-iii]
\begin{enumerate}
\item $t_2\geq 0$, $s_2=0$, $t_1=1$,
$(s_1+k)+k(t_2-1)\leq m-b_j(n-t_2)$ where $(k,m,n-1)\in T$
and $j=\nu(k,m,n-1)$.
\item $t_2=0$, $s_2=1$, $s_1\geq 2$ and let $p=km+t_1-n(<m)$,
$(k,p,m-1)\in T$ and $j=\nu(k,p,m-1)$.
\begin{enumerate}
\item $m\geq kp$ and $t_1+ks_1\leq p$.
\item $kp>m$ and $t_1+ks_1\leq p-b_j(m-s_1)$.
\end{enumerate}
\end{enumerate}
\item[17-iv]
\begin{enumerate}
\item $t_2=s_2=1$ and let $p=km+t_1-n$, $(k,p,m-1)\in T$ and
$j=\nu(k,p,m-1)$.
\begin{enumerate}
\item $m-1\geq kp$ and $(k+t_1-2)+k(k+s_1-2)\leq p$.
\item $kp> m$ and $(k+t_1-2)+k(k+s_1-2)\leq p-b_j(n-k-s_1)$.
\end{enumerate}
\end{enumerate}
\end{itemize}
\end{enumerate}
\end{itemize}

\BDEF\label{def_412}
Let $T$ be the set of triplets $(k,m,n)\in\NN^3$ satisfying
$k\geq 2$, $n>m\geq2$ and $k+m^2+n^2-2 > kmn$.
For $(k,m,n)\in T$ there exists a $j\in\NN$ such that
$(\GL_1^k\times\SL_m\times\SL_n,(\omega_1\otimes\omega_1)^{\oplus k})$
is transformed to a trivial prehomogeneous module by $j$ castling
transformations. This number $j$ is uniquely determined if we use
only castling transformations decreasing the module's dimension.
This unique $j$ will be denoted by $\nu(k,m,n)$. Thus
we obtain a map $\nu:T\Ra\NN$. For example, $\nu(k,m,n)=0$ if
and only if $mk\leq n$.\index{$\nu(k,m,n)$}\index{$(k,m,n)\in T$}
We define $(a_i)$ to be the sequence
\[
a_{-1} = -1,\qquad a_0 = 0,\qquad a_i = k a_{i-1} - a_{i-2}\ \mathrm{for}\ i>0.
\]
\EDEF

There are some cases of the form KII IV belonging neither to
KII-16 nor KII-17, but to KII I instead. These are the cases (4.1-i),
(4.1-ii), (4.7) and (4.8) from Section  4.2 in Kimura et al.
\cite{kimuraII}. We will list them here for the sake of completeness.
\begin{enumerate}
\item %
$(\GL_1^{1+s+t}\times\SL_m\times\SL_n,
(\omega_1\otimes 1)^{(*)\oplus s}\oplus(\omega_1\otimes\omega_1)\oplus(1\otimes(\omega_1^*\oplus\omega_1^{\oplus t-1})))$\\
with $k=1$, $t\geq 1$, $n\geq m+t-1$ and $s\leq m$. This is the case KII I-3.
\item %
$(\GL_1^{k+s+t}\times\SL_m\times\SL_n,
((\omega_1^{\oplus s-1}\oplus\omega_1^{(*)})\otimes 1)\oplus(\omega_1\otimes\omega_1)^{\oplus k}\oplus(1\otimes(\omega_1^*\oplus\omega_1^{\oplus t-1})))$\\
with $k\geq 2$, $t\geq 1$, $n\geq km+t-1$ and $s+k\leq m+1$. This is the case KII I-3.
\item %
$(\GL_1^{k+s+t}\times\SL_m\times\SL_n,
(\omega_1^{(*)}\otimes 1)^{\oplus s}\oplus(\omega_1\otimes\omega_1)^{\oplus k}\oplus(1\otimes\omega_1)^{\oplus t})$\\
with $n\geq km+t$, and $(\GL_1^s\times\SL_m,\omega_1^{(*)\oplus s})$
is a simple prehomogeneous module.
This is the case KII I-1.
\item %
$(\GL_1^{k+s+t}\times\SL_m\times\SL_n,
(\omega_1^{(*)}\otimes 1)^{\oplus s}\oplus(\omega_1\otimes\omega_1)^{\oplus k}\oplus(1\otimes\omega_1)^{\oplus t})$\\
with $t\geq 3$, $n=km+t-1$, and
$(\GL_1^{k+s}\times\SL_m,\omega_1^{*\oplus k}\oplus\omega_1^{(*)\oplus s})$ is a simple
pre\-homogeneous module.
This is the case KII I-2.
\end{enumerate}

\section{Tables of \'etale modules}\label{appendix_special}

In this chapter, we present the \'etale modules from part IV
of Globke \cite{globke}.
Note that this is not claimed to be a complete classification.

\subsection{\'Etale modules with center $\GL_1$}
% \'Etale with torus GL_1
\subsubsection{One-simple \'etale modules with center $\GL_1$}
\begin{itemize}
\item SK I-4: $(\GL_2,3\omega_1,\Sym^3\CC^2)$.
\item Ks I-2: $(\GL_1\times\SL_n,\mu\otimes\omega_1^{\oplus n},(\CC^n)^{\oplus n})$.
\end{itemize}
\subsubsection{Two-simple \'etale modules with center $\GL_1$}
\begin{itemize}
\item SK I-8: $(\SL_3\times\GL_2,2\omega_1\otimes\omega_1,\Sym^2\CC^3\otimes\CC^2)$.
\item SK I-11: $(\SL_5\times\GL_4,\omega_2\otimes\omega_1,\bigwedge^2\CC^5\otimes\CC^4)$.
\item KII I-2: $(\GL_1\times G\times\SL_n,
(\rho_1\otimes\omega_1)\oplus\ldots\oplus(\rho_k\otimes\omega_1)
\oplus (\rho_{k+1}^*\otimes 1)\oplus\ldots\oplus(\rho_l^*\otimes 1))$,\\
with $n=-1+\sum_{i=1}^k \dim(\rho_i)$ and
$(\GL_1\times G, \rho_1\oplus\ldots\oplus \rho_l)$ an \'etale
module for a simple group $G$.
\item KII IV-16-iv (a): $(\GL_m\times\SL_{m+1}, (\omega_1\otimes\omega_1)\oplus(1\otimes\omega_1^*)^{\oplus m+1})$, $m\geq 2$.
\end{itemize}

\subsection{\'Etale modules for groups with factor $\Sp_m$ }
% \'Etale with Sp_n
\begin{itemize}
\item KI I-16: $(\GL_1^2\times\Sp_2\times\SL_3,(\omega_1\otimes\omega_1)\oplus(\omega_2\otimes 1)\oplus(1\otimes\omega_1^*),(\CC^4\otimes\CC^3)\oplus V^5\oplus\CC^3)$.
\item KI I-18: $(\GL_1^3\times\Sp_2\times\SL_2,(\omega_2\otimes\omega_1)\oplus(\omega_1\otimes 1)\oplus(1\otimes\omega_1),(V^5\otimes\CC^2)\oplus \CC^4\oplus\CC^2)$.
\item KI I-19: $(\GL_1^3\times\Sp_2\times\SL_4,(\omega_2\otimes\omega_1)\oplus(\omega_1\otimes 1)\oplus(1\otimes\omega_1^*),(V^5\otimes\CC^4)\oplus \CC^4\oplus\CC^4)$.
\end{itemize}

% all \'Etale modules
\subsection{One-simple \'etale modules}
\begin{itemize}
% SK
\item SK I-4: $(\GL_2,3\omega_1,\Sym^3\CC^2)$.
% Ks
\item Ks I-2: $(\GL_1\times\SL_n,\mu\otimes\omega_1^{\oplus n},(\CC^n)^{\oplus n})$.
\item Ks I-3: $(\GL_1^{n+1}\times\SL_n,\omega_1^{\oplus n+1},(\CC^n)^{\oplus n+1})$.
\item Ks I-4: $(\GL_1^{n+1}\times\SL_n,\omega_1^{\oplus n}\oplus\omega_1^*,(\CC^n)^{\oplus n}\oplus\CC^{n*})$.
\item Ks I-11 for $n=2$: $(\GL_1^2\times\SL_2,2\omega_1\oplus
\omega_1,\Sym^2\CC^2\otimes\CC^2)$.
\end{itemize}

\subsection{Two-simple \'etale modules}
\begin{itemize}
% SK
\item SK I-8: $(\SL_3\times\GL_2,2\omega_1\otimes\omega_1,\Sym^2\CC^3\otimes\CC^2)$.
\item SK I-11: $(\SL_5\times\GL_4,\omega_2\otimes\omega_1,\bigwedge^2\CC^5\otimes\CC^4)$.
% KI
\item KI I-1: $(\GL_1^2\times\SL_4\times\SL_2,(\omega_2\otimes\omega_1)\oplus(\omega_1\otimes\omega_1),(\bigwedge^2\CC^4\otimes\CC^2)\oplus(\CC^4\otimes\CC^2))$.
\item KI I-2: $(\GL_1^2\times\SL_4\times\SL_2,(\omega_2\otimes\omega_1)\oplus(\omega_1\otimes1)\oplus(\omega_1\otimes1),(\bigwedge^2\CC^4\otimes\CC^2)\oplus\CC^4\oplus\CC^4)$.
\item KI I-6: $(\GL_1^3\times\SL_5\times\SL_2,(\omega_2\otimes\omega_1)\oplus(\omega_1^*\otimes 1)\oplus(\omega_1^{(*)}\otimes 1),(\bigwedge^2\CC^5\otimes\CC^2)\oplus\CC^{5*}\oplus\CC^{5(*)})$.
\item KI I-16: $(\GL_1^2\times\Sp_2\times\SL_3,(\omega_1\otimes\omega_1)\oplus(\omega_2\otimes 1)\oplus(1\otimes\omega_1^*),(\CC^4\otimes\CC^3)\oplus V^5\oplus\CC^3)$.
\item KI I-18: $(\GL_1^3\times\Sp_2\times\SL_2,(\omega_1\otimes\omega_1)\oplus(\omega_1\otimes 1)\oplus(1\otimes\omega_1),(\CC^4\otimes\CC^2)\oplus \CC^4\oplus\CC^2)$.
\item KI I-19: $(\GL_1^3\times\Sp_2\times\SL_4,(\omega_1\otimes\omega_1)\oplus(\omega_1\otimes 1)\oplus(1\otimes\omega_1^*),(\CC^4\otimes\CC^4)\oplus \CC^4\oplus\CC^4)$.
% KII I
\item KII I-1: $(\GL_1^j\times G\times \GL_n,
((\sigma_1\oplus\ldots\oplus\sigma_s)\otimes\omega_1)
\oplus ((\rho_1\oplus\ldots\oplus\rho_l)\otimes 1))$,\\
with $n=\sum_{i=1}^s\dim(\rho_i)$ and
$(\GL_1^j\times G,\rho_1\oplus\ldots\oplus\rho_l)$ an \'etale module
for a simple group $G$, $1\leq j\leq l$.
\item KII I-2: $(\GL_1^{j+t}\times G\times\SL_n,
((\rho_1\oplus\ldots\oplus\rho_k)\otimes\omega_1)
\oplus((\rho_{k+1}^*\oplus\ldots\oplus\rho_l^*)\otimes 1)
\oplus(1\otimes\omega_1^{\oplus t}))$,\\
with $n=t-1+\sum_{i=1}^k\dim(\rho_i)$, $1\leq j\leq l$,
and $(\GL_1^j\times G,\rho_1\oplus\ldots\oplus\rho_l)$ an \'etale
module for a simple group $G$.
\item KII I-3: $(\GL_1^{j+t}\times G\times\SL_n,
((\rho_1\oplus\ldots\oplus\rho_k)\otimes\omega_1)
\oplus((\rho_{k+1}\oplus\ldots\oplus\rho_l)\otimes 1)
\oplus(1\otimes(\omega_1^{\oplus t-1}\oplus\omega_1^*)))$,\\
with $n=t-1+\sum_{i=1}^k\dim(\rho_i)$, $1\leq j\leq l$, and
$(\GL_1^j\times G,\rho_1\oplus\ldots\oplus\rho_l)$ an \'etale
module for a simple group $G$.
% KII II
\item KII II-4-i (b):
$(\GL_1^3\times\SL_2\times\SL_3,(\omega_1\otimes\omega_1)\oplus(1\otimes2\omega_1^{(*)})\oplus(\omega_1^{(*)}\otimes 1))$.
\item KII II-4-ii (a):
$(\GL_1^5\times\SL_3\times\SL_4,(\omega_1\otimes\omega_1)\oplus(1\otimes\omega_2^{(*)})\oplus((\omega_1\oplus\omega_1)^{(*)}\otimes 1)\oplus(1\otimes\omega_1^{(*)}))$.
\item KII II-4-iii (d):
$(\GL_1^5\times\SL_2\times\SL_3,(\omega_1\otimes\omega_1)\oplus(1\otimes\omega_2)\oplus((\omega_1\oplus\omega_1)^{(*)}\otimes1)\oplus(1\otimes\omega_1))$.
\item KII II-4-iii (f):
$(\GL_1^5\times\SL_2\times\SL_3,(\omega_1\otimes\omega_1)\oplus(1\otimes\omega_2)\oplus(\omega_1\otimes1)\oplus(\omega_1^*\otimes1)\oplus(1\otimes\omega_1))$.
\item KII II-4-iii (g):
$(\GL_1^4\times\SL_2\times\SL_3,(\omega_1\otimes\omega_1)\oplus(1\otimes\omega_2)\oplus(1\otimes\omega_1)\oplus(1\otimes\omega_1^*))$.
\item KII II-4-iii (h):
$(\GL_1^4\times\SL_2\times\SL_3,(\omega_1\otimes\omega_1)\oplus(1\otimes\omega_2)\oplus(1\otimes\omega_1^*)\oplus(1\otimes\omega_1^*))$.
\item KII II-4-iii (n):
$(\GL_1^5\times\SL_2\times\SL_3,(\omega_1\otimes\omega_1)\oplus(1\otimes\omega_2^*)\oplus((\omega_1\oplus\omega_1)^{(*)}\otimes1)\oplus(1\otimes\omega_1^{(*)}))$.
\item KII II-4-iii (p):
$(\GL_1^4\times\SL_2\times\SL_3,(\omega_1\otimes\omega_1)\oplus(1\otimes\omega_2^*)\oplus(1\otimes\omega_1^*)\oplus(1\otimes\omega_1^*))$.
\item KII II-4-iii (q):
$(\GL_1^5\times\SL_2\times\SL_3,(\omega_1\otimes\omega_1)\oplus(1\otimes\omega_2^*)\oplus(\omega_1\otimes1)\oplus(\omega_1^*\otimes1)\oplus(1\otimes\omega_1^{(*)}))$.
\item KII II-5-i (b):
$(\GL_1^3\times\SL_2\times\SL_3,(2\omega_1\otimes\omega_1)\oplus(1\otimes\omega_2^{(*)})\oplus(\omega_1\otimes1))$.
\item KII II-5-ii (b):
$(\GL_1^4\times\SL_2\times\SL_3,(\omega_1\otimes\omega_1)\oplus(1\otimes\omega_2)\oplus(2\omega_1\otimes1)\oplus(1\otimes\omega_1))$.
\item KII II-5-iii (e):
$(\GL_1^4\times\SL_2\times\SL_3,(\omega_1\otimes\omega_1)\oplus(1\otimes\omega_2^*)\oplus(2\omega_1\otimes1)\oplus(1\otimes\omega_1^{(*)}))$.
\item KII II-5-iv (a):
$(\GL_1^3\times\SL_2\times\SL_3,(\omega_1\otimes\omega_1)\oplus(1\otimes\omega_2^*)\oplus(1\otimes\omega_2^*))$.
\item KII II-6 (b):
$(\GL_1^3\times\SL_3\times\SL_5,(\omega_1\otimes\omega_1)\oplus(1\otimes\omega_2)\oplus(1\otimes\omega_2))$.
\item KII II-9 (a):
$(\GL_1^4\times\SL_6\times\SL_7,(\omega_1\otimes\omega_1)\oplus(1\otimes\omega_2^*)\oplus(1\otimes\omega_1)\oplus(1\otimes\omega_1))$.
\item KII II-10: $(\GL_1^{r+1}\times\SL_n\times\SL_n,(\omega_1\otimes\omega_1)
\oplus (\rho_1\otimes1)\oplus\ldots\oplus(\rho_k\otimes1)
\oplus(1\otimes\rho_{k+1}^*)\oplus\ldots\oplus(1\otimes\rho_r^*))$,
where $(\GL_1^r\times\SL_n,\rho_1\oplus\ldots\oplus\rho_r)$ is an 
\'etale module.
% KII III
\item KII III-12:
$(\GL_1^2\times\SL_4\times\SL_8,(\omega_2\otimes\omega_1)\oplus(\omega_1\otimes\omega_1))$.
\item KII III-13-i (d):
$(\GL_1^5\times\SL_2\times\SL_3,(\omega_1\otimes\omega_1)\oplus(1\otimes\omega_1^*)\oplus(1\otimes\omega_1^*)\oplus(\omega_2^{(*)}\otimes1)\oplus(1\otimes\omega_1^*))$.
\item KII III-13-i (e):
$(\GL_1^5\times\SL_2\times\SL_3,(\omega_1\otimes\omega_1)\oplus(1\otimes\omega_1^*)\oplus(1\otimes\omega_1^*)\oplus(\omega_2^*\otimes1)\oplus(1\otimes\omega_1))$.
\item KII III-13-i (g):
$(\GL_1^5\times\SL_2\times\SL_3,(\omega_1\otimes\omega_1)\oplus(1\otimes\omega_1^*)\oplus(1\otimes\omega_1^*)\oplus(\omega_2\otimes1)\oplus(1\otimes\omega_1))$.
\item KII III-13-ii (a):
$(\GL_1^5\times\SL_2\times\SL_3,(\omega_1\otimes\omega_1)\oplus(1\otimes\omega_1)\oplus(1\otimes\omega_1)\oplus(1\otimes\omega_1)\oplus(\omega_2\otimes1))$.
\item KII III-13-ii (b):
$(\GL_1^5\times\SL_2\times\SL_3,(\omega_1\otimes\omega_1)\oplus(1\otimes\omega_1)\oplus(1\otimes\omega_1)\oplus(1\otimes\omega_1)\oplus(\omega_2^*\otimes1))$.
% KII IV 16
\item KII IV-16-i: $(\GL_1^{n+2}\times\SL_n\times\SL_n,
(\omega_1\otimes1)^{\oplus s_1}\oplus(\omega_1^*\otimes1)^{\oplus s_2}\oplus
(\omega_1\otimes\omega_1)\oplus(1\otimes\omega_1)^{\oplus t_1}\oplus(1\otimes\omega_1^*)^{\oplus t_2})$,\\
with $s_1+t_2=n$ and $s_2+t_1=1$.
\item KII IV-16-i: $(\GL_1^{n+2}\times\SL_n\times\SL_n,
(\omega_1\otimes1)^{\oplus s_1}\oplus
(\omega_1\otimes\omega_1)\oplus(1\otimes\omega_1^*)^{\oplus t_2})$,\\
with $s_1+t_2=n+1$.
\item KII IV-16-i: $(\GL_n\times\GL_n,
(\omega_1\otimes1)^{\oplus s_1}\oplus
(\omega_1\otimes\omega_1)\oplus(1\otimes\omega_1^*)^{\oplus t_2})$,\\ with $s_1+t_2=n$.
\item KII IV-16-ii (a): $(\GL_n\times\GL_n,
(\omega_1\otimes1)^{\oplus s_1}\oplus
(\omega_1\otimes\omega_1)^{\oplus k}\oplus(1\otimes\omega_1^*)^{\oplus t_2})$,\\
with $s_1+kt_2=m$.
\item KII IV-16-ii (b): $(\GL_n\times\GL_n,
(\omega_1^*\otimes1)^{\oplus s_2}\oplus
(\omega_1\otimes\omega_1)\oplus(1\otimes\omega_1)^{\oplus t_1})$,\\ with $s_2+kt_1=m$.
\item KII IV-16-iii: $(\GL_1^{t_1-1}\times\GL_m\times\GL_n,
(\omega_1^*\otimes1)^{\oplus s_2}\oplus
(\omega_1\otimes\omega_1)\oplus(1\otimes\omega_1)
\oplus(\mu\otimes1\otimes\omega_1)^{\oplus t_1-1}$,\\ with $s_2+t_1=m+1$.
\item KII IV-16-iv (a):
$(\GL_1^{m+1}\times\SL_m\times\SL_n,
(\mu\otimes\omega_1^*\otimes1)\oplus
(\omega_1\otimes\omega_1)\oplus(\mu\otimes 1\otimes\omega_1^*)^{\oplus m})$,\\
with $n=m+1$.
\item KII IV-16-iv (a):
$(\GL_m\times\SL_n,
(\omega_1\otimes\omega_1)\oplus(1\otimes\omega_1^*)^{\oplus m})$,\\
with $n=m+1$.
\item KII IV-16-iv (a):
$(\GL_1^{m+2}\times\SL_m\times\SL_n,
(\omega_1\otimes\omega_1)\oplus(1\otimes\omega_1^*)^{\oplus m+1})$,\\
with $n=m+1$.
\item KII IV-16-iv (a):
$(\GL_1^{m+1}\times\SL_m\times\SL_n,
(\mu\otimes\omega_1\otimes1)\oplus
(\omega_1\otimes\omega_1)\oplus(\mu\otimes 1\otimes\omega_1^*)^{\oplus m})$,\\
with $n=m+1$.
\item KII IV-16-iv (c): $(\GL_1^{t_2}\times\GL_m\times\GL_n,
(\omega_1\otimes1)^{\oplus s_1}\oplus
(\omega_1\otimes\omega_1)\oplus(1\otimes\omega_1)\oplus(\mu\otimes 1\otimes\omega_1^*)^{\oplus t_2})$,\\
with $s_1+t_2=m$ and $n=m+1$.
% KII IV 17
\item KII IV-17 (1a):
$(\GL_m\times\GL_n, (\omega_1\otimes1)^{\oplus s}\oplus(\omega_1\otimes\omega_1)^{\oplus k}\oplus(1\otimes\omega_1^*)^{\oplus t})$,\\
with $s a_{j+1}+t a_{j+2}=m_0$ and $n_0=km_0$.
\end{itemize}
Let $p=km+t-n$ and $q=kp-m=k^2m+kt-kn-m$.
\begin{itemize}
\item KII IV-17 (1b):
$(\GL_m\times\GL_n, (\omega_1\otimes1)^{\oplus s}\oplus(\omega_1\otimes\omega_1)^{\oplus k}\oplus(1\otimes\omega_1)^{\oplus t})$,\\
with $kp>m$, $s>0$, $t a_{j+1}+s a_{j+2}=\tilde{m}_0$ and $\tilde{n}_0=k\tilde{m}_0$
for $\tilde{m}_0=a_{j+1} p-a_j m$,
$\tilde{n}_0=a_j p-a_{j-1} m$ and $j=\nu(k,p,m)$.
\item KII IV-17 (1b):
$(\GL_m\times\GL_n, (\omega_1\otimes\omega_1)^{\oplus k}\oplus(1\otimes\omega_1)^{\oplus t})$,\\
with $kp>m$, $kq=p$, $2\leq t$ and $kt=q$.
\item KII IV-17 (1b):
$(\GL_1^{k+1}\times\SL_m\times\SL_n, (\omega_1\otimes\omega_1)^{\oplus k}\oplus(1\otimes\omega_1))$,\\
with $kp>m$, $k=q+1$ and $p=q^2+q$.
\item KII IV-17 (1b):
$(\GL_m\times\GL_n, (\omega_1\otimes\omega_1)^{\oplus k}\oplus(1\otimes\omega_1))$,\\
with $kp>m$, $k=q$ and $p=q^2$.
\item KII IV-17 (1b):
$(\GL_m\times\GL_n, (\omega_1\otimes\omega_1)^{\oplus k}\oplus(1\otimes\omega_1)^{\oplus t})$,\\
with $kp>m$, $kq>p$, $t a_{j+2}=\tilde{m}_0$
and $k\tilde{m}_0=\tilde{n}_0$ for
$\tilde{m}_0= a_{j+1} q - a_j p$,
$\tilde{n}_0 = a_j q - a_{j-1} p$
and $j=\nu(k,q,p)$.
\item KII IV-17 (1b):
$(\GL_1^{1+k+t}\times\SL_m\times\SL_n, (\omega_1\otimes1)\oplus(\omega_1\otimes\omega_1)^{\oplus k}\oplus(1\otimes\omega_1)^{\oplus t})$,\\
with $kp=m$ and $t+k=p+1$.
\item KII IV-17 (1b):
$(\GL_m\times\GL_n, (\omega_1\otimes1)\oplus(\omega_1\otimes\omega_1)^{\oplus k}\oplus(1\otimes\omega_1)^{\oplus t})$,\\
with $kp=m$ and $t+k=p$.
\item KII IV-17 (1b):
$(\GL_m\times\GL_n, (\omega_1\otimes1)^{\oplus s}\oplus(\omega_1\otimes\omega_1)^{\oplus k}\oplus(1\otimes\omega_1)^{\oplus t})$,\\
with $kp=m$, $2\leq s$ and $t+ks=p$.
\end{itemize}
Let $p=km+t-n$, $q=kp+s-m$ and $r=kq-p$.
\begin{itemize}
\item KII IV-17 (1c):
$(\GL_m\times\GL_n, (\omega_1^*\otimes1)^{\oplus s}\oplus(\omega_1\otimes\omega_1)^{\oplus k}\oplus(1\otimes\omega_1)^{\oplus t})$,\\
with $kp>m$, $kq>p$, $t>0$,
$s a_{j+1}+t a_{j+2}=\tilde{m}_0$ and $\tilde{n}_0=k\tilde{m}_0$
for $\tilde{m}_0=a_{j+1} q-a_j p$,
$\tilde{n}_0=a_j q-a_{j-1} p$ and $j=\nu(k,q,p)$.
\item KII IV-17 (1c):
$(\GL_m\times\GL_n, (\omega_1^*\otimes1)^{\oplus s}\oplus(\omega_1\otimes\omega_1)^{\oplus k})$,\\
with $kp>m$, $kq>p$, $kr=q$, $2\leq s$ and $ks=r$.
\item KII IV-17 (1c):
$(\GL_1^{k+1}\times\SL_m\times\SL_n, (\omega_1^*\otimes1)\oplus(\omega_1\otimes\omega_1)^{\oplus k})$,\\
with $kp>m$, $kq>p$, $k=r+1$ and $q=r^2+r$.
\item KII IV-17 (1c):
$(\GL_m\times\GL_n, (\omega_1^*\otimes1)\oplus(\omega_1\otimes\omega_1)^{\oplus k})$,\\
with $kp>m$, $kq>p$, $k=r$ and $q=r^2$.
\item KII IV-17 (1c):
$(\GL_m\times\GL_n, (\omega_1^*\otimes1)^{\oplus s}\oplus(\omega_1\otimes\omega_1)^{\oplus k})$,\\
with $kp>m$, $kq>p$, $kr>q$, $s a_{j+2}=\tilde{m}_0$
and $k\tilde{m}_0=\tilde{n}_0$ for
$\tilde{m}_0= a_{j+1} r - a_j q$,
$\tilde{n}_0 = a_j r - a_{j-1} q$
and $j=\nu(k,r,q)$.
\item KII IV-17 (1c):
$(\GL_1^{s+k+1}\times\SL_m\times\SL_n, (\omega_1^*\otimes1)^{\oplus s}\oplus(\omega_1\otimes\omega_1)^{\oplus k}\oplus(1\otimes\omega_1))$,\\
with $kp>m$, $kq=p$ and $s+k=q+1$.
\item KII IV-17 (1c):
$(\GL_m\times\GL_n, (\omega_1^*\otimes1)^{\oplus s}\oplus(\omega_1\otimes\omega_1)^{\oplus k}\oplus(1\otimes\omega_1)$,\\
with $kp>m$, $kq=p$ and $s+k=q$.
\item KII IV-17 (1c):
$(\GL_m\times\GL_n, (\omega_1^*\otimes1)^{\oplus s}\oplus(\omega_1\otimes\omega_1)^{\oplus k}\oplus(1\otimes\omega_1)^{\oplus t})$,\\
with $kp>m$, $kq=p$, $2\leq t$ and $s+kt=q$.
\item KII IV-17 (1c):
$(\GL_1^{s}\times\SL_m\times\GL_n, (\omega_1^*\otimes1)^{\oplus s}\oplus(\omega_1\otimes\omega_1)^{\oplus k})$,\\
with $k=p$ and $kp+s-1=m$.
\item KII IV-17 (1c):
$(\GL_1^{s+k}\times\SL_m\times\SL_n, (\omega_1^*\otimes1)^{\oplus s}\oplus(\omega_1\otimes\omega_1)^{\oplus k})$,\\
with $k=p+1$ and $kp+s-1=m$.
\item KII IV-17 (1c):
$(\GL_1^{s+k+1}\times\SL_m\times\SL_n, (\omega_1^*\otimes1)^{\oplus s}\oplus(\omega_1\otimes\omega_1)^{\oplus k}\oplus(1\otimes\omega_1))$,\\
with $k=p$ and $kp+s-1=m$.
\item KII IV-17 (1c):
$(\GL_m\times\GL_n, (\omega_1^*\otimes1)^{\oplus s}\oplus(\omega_1\otimes\omega_1)^{\oplus k})$,\\
with $ks=p$ and $kp=m$.
\end{itemize}
Let $p=km+t-n-1$.
\begin{itemize}
\item KII IV-17 (2a):
$(\GL_m\times\GL_n, (\omega_1\otimes1)^{\oplus s}\oplus(\omega_1\otimes\omega_1)^{\oplus k}\oplus(1\otimes(\omega_1^*\oplus\omega_1^{\oplus t-1}))$,\\
with $kp>m$, $(t-2)a_{j+1}+(k+s)a_{j+2}=\tilde{m}_0$ and $\tilde{n}_0=k\tilde{m}_0$
for $\tilde{m}_0=a_{j+1} p-a_j m$,
$\tilde{n}_0=a_j p-a_{j-1} m$ and $j=\nu(k,p,m)$.
\item KII IV-17 (2a):
$(\GL_m\times\GL_n, (\omega_1\otimes1)^{\oplus s}\oplus(\omega_1\otimes\omega_1)^{\oplus k}\oplus(1\otimes(\omega_1^*\oplus\omega_1^{\oplus t-1}))$,\\
with $kp=m$ and $t-2+k(k+s)=p$.
\end{itemize}
Let $p=km+t-n$ and $q=kp+t-m-1$.
\begin{itemize}
\item KII IV-17 (2b):
$(\GL_m\times\GL_n, ((\omega_1\oplus\omega_1^{*\oplus s-1})\otimes1)\oplus(\omega_1\otimes\omega_1)^{\oplus k}\oplus(1\otimes\omega_1^{\oplus t})$,\\
with $kp>m$, $kq>p$, $(s-2)a_{j+1}+(k+t)a_{j+2}=\tilde{m}_0$ and
$\tilde{n}_0=k\tilde{m}_0$
for $\tilde{m}_0=a_{j+1} q-a_j p$,
$\tilde{n}_0=a_j q-a_{j-1} p$ and $j=\nu(k,q,p)$.
\item KII IV-17 (2b):
$(\GL_m\times\GL_n, ((\omega_1\oplus\omega_1^{*\oplus s-1})\otimes1)\oplus(\omega_1\otimes\omega_1)^{\oplus k}\oplus(1\otimes\omega_1^{\oplus t})$,\\
with $kp>m$, $kq=p$ and $s-2+k(k+t)=q$.
\item KII IV-17 (2b):
$(\GL_1^{k+s+t}\times\SL_m\times\SL_n, ((\omega_1\oplus\omega_1^{*\oplus s-1})\otimes1)\oplus(\omega_1\otimes\omega_1)^{\oplus k}\oplus(1\otimes\omega_1^{\oplus t})$,\\
with $m=s-1+kp$ and $k+t=p+1$.
\item KII IV-17 (2b):
$(\GL_1^s\times\SL_m\times\GL_n, ((\omega_1\oplus\omega_1^{*\oplus s-1})\otimes1)\oplus(\omega_1\otimes\omega_1)^{\oplus k}\oplus(1\otimes\omega_1^{\oplus t})$,\\
with $m=s-1+kp$ and $k+t=p$.
\item KII IV-17 (3a):
$(\GL_m\times\GL_n, (\omega_1\otimes1)^{\oplus s}\oplus(\omega_1\otimes\omega_1)^{\oplus k}\oplus(1\otimes(\omega_1\oplus\omega_1^{*\oplus t-1}))$,\\
with $(k+s) a_{j+1}+(t-2) a_{j+2}=\tilde{m}_0$ and
$\tilde{n}_0=k\tilde{m}_0$,
for $\tilde{m}_0=a_{j+1} m-a_j (n-1)$,
$\tilde{n}_0=a_j m-a_{j-1} (n-1)$ and $j=\nu(k,m,n-1)$.
\end{itemize}
Let $p=km+t-n$.
\begin{itemize}
\item KII IV-17 (3b)
$(\GL_m\times\GL_n, (\omega_1\otimes1)^{\oplus s}\oplus(\omega_1\otimes\omega_1)^{\oplus k}\oplus(1\otimes(\omega_1\oplus\omega_1^{*\oplus t-1}))$,\\
with $kp>m$, $(k+t) a_{j+1}+(s-2) a_{j+2}=\tilde{m}_0$ and
$\tilde{n}_0=k\tilde{m}_0$,
for $\tilde{m}_0=a_{j+1} p-a_j (m-1)$,
$\tilde{n}_0=a_j p-a_{j-1} (m-1)$ and $j=\nu(k,p,m-1)$.
\end{itemize}
Let $p=km+t-n-1$.
\begin{itemize}
\item KII IV-17 (4):
$(\GL_m\times\GL_n, ((\omega_1^{\oplus s-1}\oplus\omega_1^*)\otimes 1)\oplus(\omega_1\otimes\omega_1)^{\oplus k}\oplus(1\otimes(\omega_1^{\oplus t-1}\oplus\omega_1^*)))$,\\
with $kp>m-1$, $(k+t-2) a_{j+1}+(k+s-2) a_{j+2}=\tilde{m}_0$ and
$\tilde{n}_0=k\tilde{m}_0$
for $\tilde{m}_0=a_{j+1} p-a_j (m-1)$,
$\tilde{n}_0=a_j p-a_{j-1} (m-1)$ and $j=\nu(k,p,m-1)$.
\item KII IV-17 (4):
$(\GL_m\times\GL_n, ((\omega_1^{\oplus s-1}\oplus\omega_1^*)\otimes 1)\oplus(\omega_1\otimes\omega_1)^{\oplus k}\oplus(1\otimes(\omega_1^{\oplus t-1}\oplus\omega_1^*)))$,\\
with $m-1=kp$ and $(k+t-2)+k(k+s-2)=p$.
\end{itemize}

% === Literature ===

\end{document}